\documentclass[11pt]{amsart}
\usepackage[left=0.9in,right=0.9in,top=0.7in,bottom=0.7in,includeheadfoot]{geometry}
\usepackage{mathtools,amssymb,amsthm}
\usepackage{mathpazo}
\usepackage{mathrsfs}
\usepackage{microtype}
\usepackage{indentfirst}
\usepackage{graphicx}
\usepackage{xcolor}
\definecolor{citeblue}{RGB}{0,0,139}
\usepackage[colorlinks=true,linkcolor=citeblue,urlcolor=citeblue,
            citecolor=citeblue,anchorcolor=citeblue]{hyperref}
\hypersetup{
  pdftitle={Boundary branch points of Douglas--Rado disks in three dimensions},
  pdfsubject={Boundary regularity of Douglas--Rado minimizers}
}

\theoremstyle{plain}
\newtheorem{theorem}{Theorem}[section]
\numberwithin{equation}{section}
\newtheorem{lemma}[theorem]{Lemma}

\newtheorem*{openproblem}{Open Problem}
\theoremstyle{remark}
\newtheorem*{remark}{Remark}
\newtheorem{example}{Example}[section]

\newcommand{\D}{\mathbb{D}}
\newcommand{\HH}{\mathbb{H}}
\newcommand{\R}{\mathbb{R}}
\newcommand{\C}{\mathbb{C}}
\newcommand{\clD}{\overline{\D}}
\newcommand{\A}{\mathscr{A}}
\newcommand{\eps}{\varepsilon}
\newcommand{\inner}[1]{\left\langle#1\right\rangle}
\DeclareMathOperator{\diver}{div}

\title{Boundary branch points of Douglas--Rad\'o disks in $\R^3$}
\author[Gao]{Rui Gao}
\address{School of Mathematics and Statistics, Shaanxi Normal University, 620 West Chang'an Avenue, Xi'an, Shaanxi, 710119, P. R. China}
\email{rui.gao@snnu.edu.cn}
\date{\today}
\begin{document}
\begin{abstract}
We prove that every Douglas--Rad\'o minimizer spanning a $C^{1,\alpha}$ Jordan curve in $\R^3$, with $0<\alpha<1$, is a $C^{1,\alpha}$ conformal immersion up to the boundary. We also construct a planar Douglas--Rad\'o minimizer that spans a $C^1$ Jordan curve and has a boundary branch point.
\end{abstract}
\maketitle
% \tableofcontents

\section{Introduction}\label{sec:introduction}

The classical Plateau problem asks for the existence of a surface of least parametric area spanning a given Jordan curve in $\R^3$. For surfaces of disk type, this problem was resolved independently in the seminal works of Douglas~\cite{douglas} and Rad\'o~\cite{rado}. For a rectifiable Jordan curve, a Douglas--Rad\'o minimizer is realized as a conformal harmonic map from the unit disk that minimizes parametric area among all admissible disk-type surfaces spanning the given Jordan curve. The regularity of the resulting disk naturally involves both the differentiability of this parametrization and the rank of its differential. While harmonicity immediately ensures interior smoothness, boundary regularity depends on the smoothness of the prescribed curve. In particular, for a regular $C^{1,\alpha}$ Jordan curve with $0<\alpha<1$, the conformal parametrization is $C^{1,\alpha}$ up to the boundary by
Kinderlehrer~\cite{kind}, see also Hildebrandt~\cite[\S3.2, Theorem~1]{hildreg} and references therein. At points of the closed disk where the differential exists, the conformal parametrization fails to be an immersion exactly when the differential vanishes, such points are called \textit{branch points}.

For least-area disks in $\R^3$, Osserman~\cite{osserman} ruled out true interior branch points, and Gulliver~\cite{gull} and Alt~\cite{alt1,alt2} ruled out false ones.  At the boundary, the regularity of the prescribed curve and the ambient codimension both enter the question.
Gulliver--Lesley~\cite{gl} ruled out boundary branch points along real analytic arcs in $\R^3$, and White~\cite{white} ruled out true boundary branch points along such arcs in arbitrary codimension.

Without real analyticity, true boundary branch points occur in higher codimension even for least-area disks. As noted in~\cite[p.~295]{white}, for each integer $k\geq 1$, restricting the holomorphic map $z\mapsto(z^3,z^{3k+1})$ to the closure of a smoothly bounded Jordan domain in the upper half-plane whose boundary contains a segment of the real axis around $0$ produces a calibrated disk in $\C^2$. This disk spans a $C^k$ Jordan curve and exhibits a true boundary branch point at $0$. Ekholm--White--Wienholtz~\cite[p.~212]{eww} constructed an area-minimizing example in $\R^4$ with a $C^\infty$ boundary and a boundary branch point by projecting an example of Gulliver~\cite{gull91} in $\R^6$. De Lellis--De Philippis--Hirsch--Massaccesi~\cite[Theorem~1.8(a)]{ddhm} give an example in $\R^4$ with a smooth Jordan boundary and a genuine boundary singularity in the mass-minimizing current class. De Lellis--De Philippis--Hirsch~\cite{ddh} identify that example as a genus-zero Douglas--Rad\'o solution as well.

In codimension one, Hardt--Simon~\cite[Theorem~11.1]{hs} established $C^{1,\alpha}$ boundary regularity, $0<\alpha<1$, for the support of a mass-minimizing integral current whose boundary is an oriented embedded $C^2$ submanifold of multiplicity one. In particular, combining their theorem with existence and interior regularity for such currents yields an embedded orientable minimal surface spanning any oriented $C^2$ Jordan curve in $\R^3$, with area minimized over orientable surfaces of arbitrary genus. The construction in~\cite{gull91} gives a minimal disk in $\R^3$ that spans a smooth Jordan curve and has a true boundary branch point. The example fails to minimize mass as a current, as discussed by~\cite[p.~5, footnote~2]{ddhm}. Whether it minimizes area in the Douglas--Rad\'o class was left open~\cite[\S1]{ddh}. Ekholm--White--Wienholtz~\cite[Theorem~3.2]{eww} proved that every classical minimal surface spanning a smooth Jordan curve of total curvature at most $4\pi$ is smoothly embedded up to the boundary. Gulliver--Osserman--Royden~\cite[Remark~6.22(2) and \S6.23, Question~6]{gor} distinguished ramified and unramified boundary branch points and asked whether area minimality excludes the latter. This leaves one of the longest standing open problems in minimal surface theory regarding the existence of boundary branch points for the Douglas--Rad\'o solution in codimension one. 

\begin{openproblem}[{Nitsche~\cite[p.~318, footnote~3]{nitsche}}]
Let $\Gamma\subset\R^3$ be a smooth Jordan curve. Can a Douglas--Rad\'o minimizer spanning $\Gamma$ have a boundary branch point?
\end{openproblem}

In Theorem~\ref{thm:main}, we resolve this open problem by proving that every Douglas--Rad\'o minimizer spanning a $C^{1,\alpha}$ Jordan curve in $\R^3$ ($0<\alpha<1$) has no boundary branch points. Moreover, the $C^{1,\alpha}$ regularity assumption is sharp, as we construct a planar Douglas--Rad\'o minimizer spanning a $C^1$ Jordan curve that admits a boundary branch point, see Example \ref{ex:planar-c1}.

Let $0<\alpha<1$ and let $\Gamma\subset\R^3$ be an oriented $C^{1,\alpha}$ Jordan curve.  We write $\D\subset\C$ for the unit disk and denote by $\A(\Gamma)$ the class of maps $u\in C^0(\clD,\R^3)\cap W^{1,2}(\D,\R^3)$ whose boundary trace $u|_{\partial\D}$ is a weakly monotone parametrization of $\Gamma$, see Dierkes--Hildebrandt--Sauvigny~\cite[Chapter 4]{dhs} for more details on the Plateau problem. 
For a map $u$ on a planar domain $U$, its Dirichlet energy and parametric area are
\[
 E(u,U)=\frac12\int_U |du|^2\,dx\,dy,
 \qquad
 A(u,U)=\int_U |u_x\times u_y|\,dx\,dy.
\]
We abbreviate the integrals by $E(u)$ and $A(u)$ when $U$ is the full domain of $u$. Equipped with these notations, we are ready to state the main results of this paper.

\begin{theorem}\label{thm:main}
Let $0<\alpha<1$ and let $\Gamma\subset\R^3$ be an oriented $C^{1,\alpha}$ Jordan curve. Then, every minimizer $u$ of $E$ over $\A(\Gamma)$ is a $C^{1,\alpha}$ conformal immersion on $\clD$. 
\end{theorem}

\begin{remark}
    In Example~\ref{ex:planar-c1}, we construct a planar least-area disk spanning a $C^1$ Jordan curve that exhibits a boundary branch point. In the planar setting, the Kellogg--Warschawski theorem, in the form proved by Warschawski~\cite[Theorems~1 and~2]{warschawski}, guarantees that the derivative of the Riemann map extends continuously and is nonvanishing on the closed disk whenever the boundary tangent is Dini continuous.  Because the tangent to a $C^{1,\alpha}$ Jordan curve is H\"older continuous, it is automatically Dini continuous, ruling out boundary branch points in that regime. At the $C^1$ endpoint case, however, Example~\ref{ex:planar-c1} features a continuous tangent whose Dini integral diverges at $0$. For the explicit map in this example, $f'(z)\to0$ as $z\to0$, so the derivative of the Riemann map vanishes at the preimage of $0$.
\end{remark}

\subsection{Outline of the proof}\label{sec:proof-outline}

We adapt the strategy of exchanging intersecting sheets, which is used to exclude interior branch points in~\cite{osserman} and comparing graph areas through Dirichlet energy, as in
Micallef--White~\cite[Theorem~2.1]{mw}. In the homogeneous model for a
real analytic boundary, the sheets cross transversely, see~\cite[Proposition~5]{white}. The main obstruction in proving Theorem \ref{thm:main} at a $C^{1,\alpha}$ boundary is that the graph
heights $u^3_1,u^3_2$ and their gradient estimates \eqref{eq:decay} allow two distinct graph heights 
to agree to infinite order at the vertex of their common sector. We cannot infer from the boundary expansion alone that an arc of transverse intersection
reaches the vertex.  Our idea is to use an energy comparison to decide where to join the sheets in the parameter domain. We keep the average height fixed of two sector sheets. In the isothermal coordinates obtained in Subsection~\ref{sec:transfer}, their half difference becomes a harmonic function $h$ that vanishes faster than linearly at the image $p_0$ of the vertex. When $h\not\equiv0$, we use this
decay to choose the branch value of a double cover $\Pi\colon\D\to\D$ by
Lemma~\ref{lem:scalar}. We prescribe $h\circ\Pi$ on one boundary arc between the lifts of $p_0$ and $-h\circ\Pi$ on the other. For this choice of branch value, the harmonic extension $g$ satisfies $E(g)<2E(h)$. Adding back the average height recovers both graph traces, so the reglued disk still spans
$\Gamma$ once monotonously. Since the quadratic bound in Lemma~\ref{lem:maj} is exact at
the original slopes, \eqref{eq:drop} and \eqref{eq:sewn} give the desired area comparison $A(\widetilde u)-A(u)\le E(g)-2E(h)<0$.

We now outline the main steps in the proof of Theorem~\ref{thm:main}.
Suppose that a minimizer $u$ has a boundary branch point $p\in\partial\D$.
We first use the boundary estimates in~\cite{kind} to obtain the local
expansion \eqref{eq:leading} of $u$ at $p$. Using this expansion, we choose a local
coordinate $t$ in which the projection $P$ onto the limiting tangent plane
has the form $P=t^q$, with $q\ge3$ odd. Two consecutive inverse branches
$z_1,z_2$ of $P$ over a planar sector $\Omega$ with vertex $0$ determine
closed sectors $S_1,S_2$ in the parameter disk that meet only at $p$.
On $S_j$, $u$ parametrizes the minimal graph $w\mapsto(w,u^3_j(w))$.
By \eqref{eq:decay}, its height is $O(|w|^{1+\alpha/q})$ and its slope
is $O(|w|^{\alpha/q})$ as $w\to0$.

We next seek a replacement of $u$ with the same two boundary traces over
$\partial\Omega$. Writing the heights as $u^3_+\pm u^3_-$, we keep their
average $u^3_+$ fixed and vary their half difference $u^3_-$. For
small gradients of $u_+^3$ and $u^3_{-}$, Lemma~\ref{lem:maj} constructs the matrix $B$ in
\eqref{eq:transfer-coefficient}. Subtracting the minimal graph equations on two sheets gives the weak equation $\diver(B\nabla u^3_-)=0$ in \eqref{eq:transfer-weak}.
The same $B$ gives a quadratic upper bound for the sum of the two area
densities, with equality at $\nabla u^3_-$ by \eqref{eq:tangent}. Thus
lowering the Dirichlet energy of the half difference, computed with $B$
on the two sheets, lowers their total area.

Since $\det B=1$, we can choose isothermal coordinates $T\colon\Omega\to\D$
in which $h=u^3_-\circ T^{-1}$ is harmonic. We choose the opening of
$\Omega$ close to $\pi$ and then shrink its radius to control the distortion
of $T$. The height estimate \eqref{eq:decay} then gives
$|h(\zeta)|\le C|\zeta-p_0|^\delta$ near $p_0=T(0)$, with $\delta>1$ as in
\eqref{eq:dist}. Suppose first that $u^3_-\not\equiv0$. On a double cover
$\Pi\colon\D\to\D$, we prescribe $h\circ\Pi$ on one boundary arc between
the two lifts of $p_0$ and $-h\circ\Pi$ on the other. By Lemma~\ref{lem:scalar},
the decay at $p_0$ lets us choose the branch value $b_0$ so that the harmonic
extension $g$ of these data satisfies $E(g)<2E(h)$. On the double cover
$\pi_b\colon\D\to\Omega$ branched over $b=T^{-1}(b_0)$, we pull $g$ back by
the lift $\widetilde T$ in \eqref{eq:transfer-lift} and add $u^3_+\circ\pi_b$.
This sum is the height of the comparison map $Y$ in \eqref{eq:newmap}. Its two boundary arcs
carry the original graph traces, and \eqref{eq:drop} gives $A(Y)<A(u^3_1)+A(u^3_2)$ even if $u^3_-$ vanishes
to infinite order at the vertex.

To insert $Y$ into a disk spanning $\Gamma$, we cut each $S_j$ along
$z_j([0,b])$ and glue opposite banks of the two slits as in
\eqref{eq:surgery-bank-pairing}. The reglued sectors form the covering
disk of $\pi_b$ inside a disk $M$ whose boundary is identified with
$\partial\D$. Since the graph traces match, we can use $Y$ on the
reglued sectors and $u$ on their complement to obtain an admissible
disk $\widetilde u\in\A(\Gamma)$. Combining \eqref{eq:drop} with
\eqref{eq:sewn} gives $A(\widetilde u)<A(u)$, contradicting the area
minimality of $u$ in \eqref{eq:mine}.

In the remaining case $u^3_-\equiv0$, we use the same construction with any $b\in\Omega$
and the comparison map $Y=(\pi_b,u^3_+\circ\pi_b)$. Both sheets parametrize the common
graph, so the replacement preserves area. A quasiconformal change of
parameters gives a weakly conformal disk $u_*$ with
$E(u_*)=A(\widetilde u)=A(u)$. Thus $u_*$ is another energy minimizer
by \eqref{eq:mine}. Near the interior ramification point of $\pi_b$, the
distinct points with local coordinates $\pm\eta$, $\eta\ne0$, have the
same projection $w=b+\eta^2$ and the same height $u^3_+(w)$. Reparametrization
preserves this failure of local injectivity, so $u_*$ cannot be an immersion
at the corresponding interior point. This contradicts the interior immersion
theorem in~\cite[Theorem~6.2]{gull} and completes the proof of Theorem \ref{thm:main}.

\subsection{Organization of the paper}\label{sec:organization}

In Section~\ref{sec:geometry}, we begin the proof outlined in
Subsection~\ref{sec:proof-outline} by showing that a minimizer $u$ of $E$
in $\A(\Gamma)$ also minimizes area in this class. If $u$ has a boundary
branch point, we obtain two sectors $S_1,S_2$ in the parameter disk on
which $u$ parametrizes minimal graphs over a common planar sector $\Omega$.

In Section~\ref{sec:cover-comparison}, we compare the sum of the two graph
areas with the area of a map $Y$ on a branched double cover of $\Omega$.
Subsection~\ref{sec:majorant} gives a quadratic upper bound for the sum of
the two graph area densities, with equality at the original slopes. In
Subsection~\ref{sec:scalar}, we use the superlinear boundary decay of a nonzero
harmonic function $h$ to choose a double cover of $\D$ on which reversing
the sign of its lifted trace on one boundary arc gives a harmonic extension
of energy less than $2E(h)$. We combine these results in
Subsection~\ref{sec:transfer}, where we use isothermal coordinates for the
equation satisfied by the height half difference $u^3_-$. This yields a map
$Y$ with the original graph traces and the strict area inequality
\eqref{eq:drop} when the heights differ.

In Section~\ref{sec:surgery}, we use $Y$ to construct an admissible disk
and complete the proof of Theorem~\ref{thm:main}. In
Subsection~\ref{sec:surgery-topology}, we cut and reglue the parameter disk
along one slit in each of $S_1,S_2$, show that the resulting surface is a
disk, and match the two boundary traversals of the cover to the original
graph traces. In Subsection~\ref{sec:surgery-sobolev}, we prove that the
resulting map $\widetilde u$ belongs to $\A(\Gamma)$ and has strictly
smaller area than $u$ when the graph heights differ. In
Subsection~\ref{sec:surgery-coincident}, we rule out coincident heights
by showing that the area-preserving replacement would give, after
reparametrization, an energy minimizer that fails to be locally injective
at an interior point. We conclude with
Example~\ref{ex:planar-c1}, in which a planar least-area disk spans a
regular $C^1$ Jordan curve and has a boundary branch point.

\medskip
\noindent\textbf{Disclosure of AI assistance.} We chose the research questions studied in this paper, formulated our theorems and lemmas, and worked out the main ideas and steps of the proofs. An AI tool helped us redraft and polish passages of the manuscript, check and correct
technical details, and prepare the TikZ code for all ten figures. We wrote the manuscript and take full responsibility for its contents, including all mathematical arguments and conclusions.

\vspace{0.5cm}
\section{Local structure at a boundary branch point}
\label{sec:geometry}

In this section, we show that if a minimizer in Theorem~\ref{thm:main} has a boundary branch point, then two sectors of the parameter disk near that point map to minimal graphs over a common planar sector.

Let $u$ be a minimizer as in Theorem~\ref{thm:main}. By the classical Plateau theory~\cite{douglas,rado}, $u$ is harmonic and weakly
conformal in $\D$, and $u|_{\partial\D}\colon\partial\D\to\Gamma$ is a
homeomorphism, see also \cite[\S4.5, Theorems~2 and~3]{dhs}. For a rectifiable Jordan curve, the area and energy infima over $\A(\Gamma)$ agree by Hildebrandt--von der Mosel~\cite{hvm99},
see also \cite[\S4.10, Theorem~1]{dhs}. Hence, we have
\begin{equation}
\label{eq:mine}
 E(u)=A(u)=\inf_{v\in\A(\Gamma)} A(v)=\inf_{v\in\A(\Gamma)} E(v).
\end{equation} 
The homeomorphic trace parametrizes each regular $C^{1,\alpha}$ boundary arc
monotonically. Thus $u\in C^{1,\alpha}(\clD,\R^3)$ by~\cite[Theorem~4]{kind}.

We write
\[
 D_r^+=\{z\in\C:|z|<r,\ \operatorname{Im}z>0\}
\]
for the upper half disk of radius $r$. Suppose $du(p)=0$ for some $p\in\partial\D$. Choose a conformal coordinate $z=x+\sqrt{-1}y$ in which $p$ corresponds to $0$ and a neighborhood of $p$ is $D_\rho^+$ together with its diameter. After an orientation preserving rigid motion of $\R^3$, we may assume that $u(p)=0$ and that the positively oriented tangent line to $\Gamma$ at $0$ coincides with the positive first axis of $\R^3$.

Let $F=(F_1,F_2,F_3)$ be the holomorphic map on $D_\rho^+$ with
$\operatorname{Re}F=u$, normalized by $F(0)=0$. Near $0$, the curve
$\Gamma$ is a $C^{1,\alpha}$ graph over the first axis in $\R^3$, and the
second and third coordinate functions have zero derivative at $0$. The oriented boundary homeomorphism $u|_{\partial \mathbb{D}} \to \Gamma$ makes $x\mapsto u^1(x)$ strictly increasing on a real interval near $0$.

By~\cite[Theorem~1 and Corollary~1]{kind}, after shrinking
the half disk $D^+_r$ there is a positive odd integer $q$ and a univalent branch
$\zeta(z)=F_1(z)^{1/q}$ with $\zeta(0)=0$. Set
$f=F\circ\zeta^{-1}$, we have
\[
 u(z)=\operatorname{Re}f(\zeta(z)),\qquad f_1(\zeta)=\zeta^q.
\]
By~\cite[Theorem~3]{kind}, the image of the real
boundary interval under $\zeta$ is a $C^{1,\alpha}$ arc near $0$ and
$H(\zeta):=\zeta^{1-q}f'(\zeta)$ is $C^{0,\alpha}$ up to that arc. The
local boundary regularity theorem for conformal maps used in~\cite[p.~738]{kind}
gives $C^{1,\alpha}$ extensions of $\zeta$ and $\zeta^{-1}$ to the respective
boundary arcs. Their
derivatives at $0$ are reciprocal.
Writing $a=\zeta'(0)\ne0$, we obtain
\[
 \zeta(z)=az+O(|z|^{1+\alpha}),\qquad
 \zeta'(z)=a+O(|z|^\alpha).
\]
Since $f_1(\zeta)=\zeta^q$, the first component of $H$ is identically $q$,
so $H(0)\ne0$. The expansion $H(\zeta(z))=H(0)+O(|z|^\alpha)$ gives
\[
 \begin{aligned}
  2u_z(z)=F'(z)&=\zeta(z)^{q-1}\zeta'(z)H(\zeta(z)) =c z^{q-1}+O(|z|^{q-1+\alpha}),
               \qquad c=a^q H(0)\ne0.
 \end{aligned}
\]
The expansion of $2u_z$, with $c\ne0$, and the condition $du(0)=0$
imply $q>1$. Since $q$ is odd, $q\ge3$. Dividing the conformality equation
$u_z\cdot u_z=0$ by $z^{2q-2}$ and passing to the limit gives $c\cdot c=0$.
Thus $\operatorname{Re}c$ and $-\operatorname{Im}c$ are orthogonal and have
the same positive length. The identity $du=\operatorname{Re}(2u_z\,dz)$ gives
the differential of $q^{-1}\operatorname{Re}(c z^q)$ as the leading term
of $du$. Integrating along radial segments and using $u(0)=0$ gives the
corresponding expansion of $u$. Rotate and scale the target so that
$q^{-1}\operatorname{Re}c$ and $-q^{-1}\operatorname{Im}c$ become the first
two coordinate vectors. In these coordinates, we have 
\begin{equation}
\label{eq:leading}
 \begin{aligned}
 u(z) =(\operatorname{Re}z^q,\operatorname{Im}z^q,0)+O(|z|^{q+\alpha}), \qquad du(z)=d(\operatorname{Re}z^q,\operatorname{Im}z^q,0)+O(|z|^{q-1+\alpha}).
 \end{aligned}
\end{equation}

Write $P=u^1+\sqrt{-1}u^2$ for the projection onto the limiting tangent plane.
By \eqref{eq:leading}, we have 
\[
 \frac{P(z)}{z^q}=1+O(|z|^\alpha),\qquad
 D\!\left(\frac{P(z)}{z^q}\right)=O(|z|^{\alpha-1}).
\]
Since $P(z)/z^q\to1$ as $z\to0$, we use the holomorphic branch of $w^{1/q}$ near $w=1$ with value $1$ there, and define
\begin{equation*}
    t(z)=z\left(\frac{P(z)}{z^q}\right)^{1/q}\quad(z\ne0), \qquad t(0)=0.
\end{equation*}
Then $P=t^q$, and the local expansions of $P(z)/z^q$ and $D(P(z)/z^q)$ give
\[
 t(z)=z+O(|z|^{1+\alpha}),\qquad Dt(z)=I+O(|z|^\alpha).
\]
Shrink $\rho$ so that $\|Dt-I\|_{\mathrm{op}}\le 1/2$.
Integrating along segments in the closed half disk gives
\[
 \tfrac12|z-z'|\le |t(z)-t(z')|\le \tfrac32|z-z'|.
\]
Thus $t$ is bi-Lipschitz on $\overline{D_\rho^+}$. Since $u$ is harmonic,
$P=u^1+\sqrt{-1}u^2$ is smooth on $D_\rho^+$, and the defining formula for $t$
gives $t\in C^\infty(D_\rho^+)$. The bound $\|Dt-I\|_{\mathrm{op}}\le 1/2$
implies $\det Dt>0$ on $D_\rho^+$. Since $t$ is injective, $t\colon D_\rho^+\to t(D_\rho^+)$ is a smooth,
orientation preserving diffeomorphism.

For every closed angular interval contained in $(0,\pi)$, the open image $t(D_\rho^+)$ contains all points $\tau$ with argument in that interval
and sufficiently small positive modulus. To see this, choose $\eta>0$
so that the closed ball of radius $\eta|\tau|$ centered at $\tau$ lies in
$D_\rho^+$ for all such $\tau$. On that ball, $|z|\le(1+\eta)|\tau|$, so
$|t(z)-z|\le C|\tau|^{1+\alpha}$ and $\|Dt(z)-I\|_{\mathrm{op}}\le C|\tau|^\alpha$.
For small $|\tau|$, the map $z\mapsto\tau+z-t(z)$ is a contraction of the
ball into itself. Its fixed point satisfies $t(z)=\tau$.

Let $\beta=1+\alpha/q$. Choose $\pi/\beta<\theta<\pi$ and
$0<\theta_0<\pi-\theta$, and for $r>0$ set
\[
 \Omega=\{w\in\C\mid 0<|w|<r,\ \theta_0<\arg w<\theta_0+\theta\}.
\]
For $w\in\overline\Omega\setminus\{0\}$, take the two $q$th roots of $w$
with arguments $\arg w/q$ and $(2\pi+\arg w)/q$. The two closed angular
intervals for these roots are disjoint and lie in $(0,\pi)$, since
$0<\theta_0<\theta_0+\theta<\pi$ and $q\ge3$. Using the inclusion proved above for $t(D_\rho^+)$, we choose $r$ small
enough that both roots lie in $t(D_\rho^+)$ for every
$w\in\overline\Omega\setminus\{0\}$. In
the $z$-coordinate, define
\[
 z_j(w)=t^{-1}\!\left(|w|^{1/q}e^{\sqrt{-1}(\arg w+2\pi(j-1))/q}\right),
 \qquad j=1,2.
\]
Set $z_j(0)=0$. The bi-Lipschitz estimate for $t$ gives continuity of $z_j$ at $0$.
Let $S_j$ be the closed sector in the parameter disk represented by
$z_j(\overline\Omega)$ in the $z$-coordinate. The sectors meet only at
$p$, and $S_j\setminus\{p\}\subset\D$. Since $P=t^q$,
$P(z_j(w))=w$, so $z_j(w)$ is the $z$-coordinate of the unique preimage
of $w$ in $S_j$. For $w\ne0$, $Dt(z_j(w))$ is invertible and
$t(z_j(w))\ne0$, so $DP(z_j(w))$ is invertible. Thus each $z_j$ extends smoothly to a neighborhood of every
$w\in\partial\Omega\setminus\{0\}$.

\begin{figure}[!htbp]
  \centering
  \includegraphics[width=.95\linewidth]{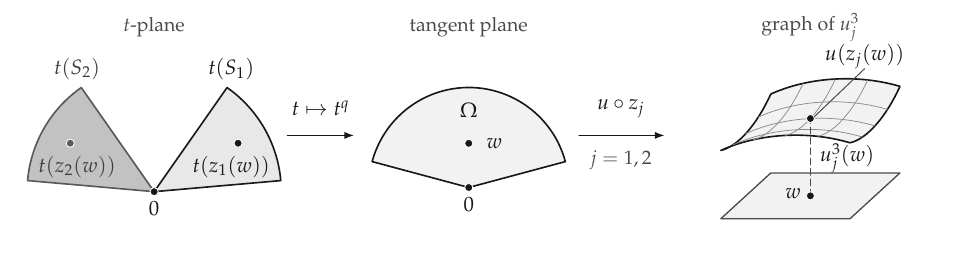}
  \caption{The two inverse branches over a fixed $w\in\Omega$, drawn for $q=3$.
  The marked points in $t(S_1)$ and $t(S_2)$ both have $q$th power $w$.
  The right panel shows the graph of $u^3_j$ near $w$, with
  $u(z_j(w))=(w,u^3_j(w))$ for $j=1,2$.}
  \label{fig:geometry-two-heights}
\end{figure}

Define the two graph heights by
\[
 u^3_j(w)=u^3(z_j(w)),\qquad w\in\overline\Omega,\quad j=1,2.
\]
Then $u(z_j(w))=(w,u^3_j(w))$ and $u^3_j(0)=0$.
Figure~\ref{fig:geometry-two-heights} shows how the two inverse branches determine the graph heights.

Over $\Omega$, the graph is minimal because $P$ is a local diffeomorphism at $z_j(w)$ and $u$ is conformal and harmonic. For $w\in\Omega$ and $z=z_j(w)$, the identity $P=t^q$ and the bound on $Dt-I$ give
\[
 \|Dz_j(w)\|_{\mathrm{op}}=\|(DP(z))^{-1}\|_{\mathrm{op}}
 \le C|t(z)|^{1-q}\le C'|z|^{1-q}.
\]
By \eqref{eq:leading}, $|u^3(z)|\le C|z|^{q+\alpha}$ and
$\|Du^3(z)\|_{\mathrm{op}}\le C|z|^{q-1+\alpha}$. Since $u^3_j=u^3\circ z_j$ and $|w|=|t(z)|^q\asymp |z|^q$, we obtain
\begin{equation}
\label{eq:decay}
 |u^3_j(w)|\le C|w|^\beta,\qquad
 |\nabla u^3_j(w)|\le C|w|^{\alpha/q},\qquad j=1,2.
\end{equation}
The gradient estimate and convexity of $\Omega$ give a Lipschitz bound of $u_j^3$
on $\Omega$ and $u^3_j$ are Lipschitz continuous on $\overline\Omega$. Write the average and half difference of the heights as
\[
 u^3_+=\frac{u^3_1+u^3_2}{2},\qquad
 u^3_-=\frac{u^3_1-u^3_2}{2}.
\]

\vspace{0.5cm}
\section{Area comparison on a branched double cover}
\label{sec:cover-comparison}

In this section, we compare the sum of the areas of the two minimal graphs over $\Omega$ from Section~\ref{sec:geometry} with the area of a
map $Y$ on a branched double cover of $\Omega$. The trace of $Y$ consists of one positive traversal of each graph boundary. If the graph heights differ, we obtain the strict area inequality \eqref{eq:drop}. If they coincide, the area of $Y$ equals the sum of
the two graph areas.

\subsection{A quadratic upper bound for the two graph areas}
\label{sec:majorant}

In this subsection, we prove Lemma~\ref{lem:maj}, which gives a quadratic
upper bound for the sum of the two graph area integrands as a function of
the half difference $v$ of the slopes, with the average slope $a$ fixed.

For $v\in\R^2$, set
\[
 W(v)=\sqrt{1+|v|^2},\qquad
 F_a(v)=W(a+v)+W(a-v).
\]

\begin{lemma}\label{lem:maj}
There exist constants $\eps_0,C>0$ such that to each pair $a,s\in\R^2$ with
$|a|,|s|\le\eps_0$ one can assign a symmetric positive definite
$2\times2$ matrix $B_{a,s}$. The assignment $(a,s)\mapsto B_{a,s}$ is Borel measurable and
satisfies
\[
 \det B_{a,s}=1,\qquad
 \|B_{a,s}-I\|_{\mathrm{op}}\le C(|a|^2+|s|^2).
\]
Moreover,
\begin{equation}\label{eq:tangent}
 \begin{aligned}
 2B_{a,s}s&=DW(a+s)-DW(a-s),\\
 F_a(v)&\le F_a(s)+\inner{B_{a,s}v,v}-\inner{B_{a,s}s,s}
              \qquad\text{for all }v\in\R^2.
 \end{aligned}
\end{equation}
\end{lemma}

\begin{proof}
We take $\eps_0\le1/32$ and use $c,C>0$ for universal constants. If $s=0$,
set $B_{a,0}=I$. Since $F_a$ is even, $DF_a(0)=0$. Also,
\[
 D^2W(v)=\frac{I}{W(v)}-\frac{v\otimes v}{W(v)^3}\le I,
\]
so $D^2F_a\le 2I$ and hence $F_a(v)\le F_a(0)+|v|^2$.

Suppose $s\ne0$. The function $F_a(v)$ depends on $v$ through $|v|^2$
and $(a\cdot v)^2$. Indeed, with $R=1+|a|^2+Q$, define
\begin{equation}\label{eq:maj-density}
 \Phi(Q,\ell)=\sqrt{R+2\sqrt{\ell}}+\sqrt{R-2\sqrt{\ell}}
             =\sqrt{2\bigl(R+\sqrt{R^2-4\ell}\bigr)}.
\end{equation}
Then $F_a(v)=\Phi(|v|^2,(a\cdot v)^2)$. If $|v|\le8$, the pair
$(|v|^2,(a\cdot v)^2)$ lies in the rectangle
$[0,64]\times[0,1/16]$. The second expression in \eqref{eq:maj-density} is smooth on a
neighborhood of this rectangle, with derivative bounds uniform for
$|a|\le1/32$.

Define
\begin{equation}\label{eq:maj-tangent-matrix}
 B_0=\Phi_Q\bigl(|s|^2,(a\cdot s)^2\bigr)I
     +\Phi_\ell\bigl(|s|^2,(a\cdot s)^2\bigr)(a\otimes a).
\end{equation}
Differentiating $F_a$ at $s$ gives $2B_0s=DF_a(s)=DW(a+s)-DW(a-s)$. At
$(|s|^2,(a\cdot s)^2)$,
\[
 \Phi_Q=\frac12\left(\frac{1}{W(a+s)}+\frac{1}{W(a-s)}\right),
 \qquad -C\le\Phi_\ell\le0.
\]
The sign of $\Phi_\ell$ follows from the second expression in
\eqref{eq:maj-density}.
For $0\le x\le4\eps_0^2$, we have $x/4\le1-(1+x)^{-1/2}\le x/2$.
Applying these bounds to $|a+s|^2$ and $|a-s|^2$ gives
\[
  \frac{|a|^2+|s|^2}{4}\le 1-\Phi_Q\le\frac{|a|^2+|s|^2}{2}.
\]
Using the bounds for $\Phi_\ell$, we obtain
\[
 \bigl(1-C(|a|^2+|s|^2)\bigr)I\le B_0
     \le\left(1-\frac{|a|^2+|s|^2}{4}\right)I.
\]
In particular, $B_0\ge I/2$ for sufficiently small $\eps_0$.

Let $e=s/|s|$ and let $P^\perp=I-e\otimes e$ be the orthogonal projection onto $s^\perp$. To make the determinant one of $B_{a,s}$ while preserving the action on
$s$, set
\begin{equation}\label{eq:maj-determinant-correction}
 \mu=\frac{1-\det B_0}{\inner{B_0e,e}},\qquad B_{a,s}=B_0+\mu P^\perp.
\end{equation}
In an orthonormal basis beginning with $e$, the correction term $\mu P^\perp$ changes the
determinant of $B_0$ by $\mu\inner{B_0e,e}$. Thus $\det B_{a,s}=1$, and
$P^\perp s=0$ gives $B_{a,s}s=B_0s$. Since the eigenvalues of $B_0$ are at
most one, $1-\det B_0$ lies between $1-\lambda_{\max}(B_0)$ and
$2\|B_0-I\|_{\mathrm{op}}$. Together with $1/2\le\inner{B_0e,e}\le1$,
this gives
\[
 \frac{|a|^2+|s|^2}{4}\le\mu\le C(|a|^2+|s|^2),\qquad
 \|B_{a,s}-I\|_{\mathrm{op}}\le C(|a|^2+|s|^2).
\]
Since $\mu\ge0$, we have $B_{a,s}\ge B_0\ge I/2$. The two eigenvalues of
$B_{a,s}$ have product one, so
\begin{equation}\label{eq:maj-ellipticity}
 \frac12 I\le B_{a,s}\le2I.
\end{equation}

For $|v|\le8$, write
$\Delta Q=|v|^2-|s|^2$ and $\Delta\ell=(a\cdot v)^2-(a\cdot s)^2$.
On the rectangle $[0,64]\times[0,1/16]$, differentiating the sum of the
two square roots in \eqref{eq:maj-density} gives
\[
 \Phi_{QQ}=-\frac14\left((R+2\sqrt{\ell})^{-3/2}
                    +(R-2\sqrt{\ell})^{-3/2}\right)\le-c<0.
\]
The other second derivatives of $\Phi$ are uniformly bounded on the same
rectangle. The segment from $(|s|^2,(a\cdot s)^2)$ to
$(|v|^2,(a\cdot v)^2)$ lies in this rectangle. By \eqref{eq:maj-tangent-matrix}, the linear term in the expansion of $\Phi$
along this segment about $(|s|^2,(a\cdot s)^2)$ is
\[
 \Phi_Q\Delta Q+\Phi_\ell\Delta\ell
   =\inner{B_0v,v}-\inner{B_0s,s}.
\]
Absorbing the mixed term into the negative term in $(\Delta Q)^2$ gives, for a smaller $c>0$,
\[
 \begin{aligned}
 F_a(v)-F_a(s)-\inner{B_0v,v}+\inner{B_0s,s} \le -c(\Delta Q)^2+C(\Delta\ell)^2.
 \end{aligned}
\]
To control $\Delta\ell$, use
$v=(v\cdot e)e+P^\perp v$ and $s=|s|e$ to expand
\[
 \begin{aligned}
 \Delta\ell-(a\cdot e)^2\Delta Q
 =2(v\cdot e)(a\cdot e)(a\cdot P^\perp v)  +(a\cdot P^\perp v)^2-(a\cdot e)^2|P^\perp v|^2.
 \end{aligned}
\]
Since $|v|\le8$, it follows that
\[
 (\Delta\ell)^2\le C|a|^4\bigl((\Delta Q)^2+|P^\perp v|^2\bigr).
\]
Substitute the bound on $(\Delta\ell)^2$ into the Taylor expansion estimate and
choose $\eps_0$ small enough that $C\eps_0^4\le c$ and
$C\eps_0^2\le1/4$. The Taylor remainder is then at most
$\frac14|a|^2|P^\perp v|^2\le\mu|P^\perp v|^2$. By
\eqref{eq:maj-determinant-correction} and $P^\perp s=0$, adding
$\mu|P^\perp v|^2$ replaces $B_0$ by $B_{a,s}$ in the quadratic bound.
This proves the inequality in \eqref{eq:tangent} for $|v|\le8$.

For $|v|>8$, we have $W(a\pm v)\le1+|a\pm v|\le 1+\eps_0+|v|$
and hence $F_a(v)\le2(1+\eps_0+|v|)$. Using $F_a(s)\ge2$,
$|s|\le\eps_0$, and \eqref{eq:maj-ellipticity}, we obtain
\[
 \begin{aligned}
 F_a(s)+\inner{B_{a,s}v,v}-\inner{B_{a,s}s,s} \ge 2-2\eps_0^2+\tfrac12|v|^2  \ge 2(1+\eps_0+|v|).
 \end{aligned}
\]
The last inequality uses $|v|>8$ and $\eps_0\le1/32$. This proves the
bound in \eqref{eq:tangent} for every $v\in\R^2$. The formula for $B_{a,s}$ is
continuous on $s\ne0$, and $B_{a,0}=I$, so the assignment $(a,s) \mapsto B_{a,s}$ is Borel
measurable. This completes the proof of Lemma \ref{lem:maj}.
\end{proof}

\subsection{Dirichlet energy comparison}\label{sec:scalar}

In this subsection, we prove the Dirichlet energy comparison in
Lemma~\ref{lem:scalar}. We choose the branch value of a double cover so that changing the sign of the lifted trace on one boundary traversal gives a harmonic extension with energy strictly below that of the lift.

\begin{lemma}\label{lem:scalar}
Let $h$ be a nonzero real-valued harmonic function on $\D$ with finite Dirichlet energy and a H\"older continuous extension to $\clD$. Suppose that for some $p_0\in\partial\D$ and $\delta>1$,
\begin{equation}\label{eq:scalar-decay}
 |h(z)|\le C|z-p_0|^\delta\qquad\text{for }z\in\clD\text{ near }p_0.
\end{equation}
There exists a holomorphic branched double cover $\pi\colon\D\to\D$ extending smoothly to $\clD$ with a single simple ramification point and branch value $b_0\in\D$. The set $\pi^{-1}(p_0)$ consists of two boundary points. On each component of $\partial\D\setminus\pi^{-1}(p_0)$, the restriction
of $\pi$ is an orientation preserving diffeomorphism onto
$\partial\D\setminus\{p_0\}$. For either ordering of these two components,
there is a harmonic function $g\in C(\clD)\cap W^{1,2}(\D)$ with boundary
values $g=h\circ\pi$ on the first component, $g=-h\circ\pi$ on the second,
and $g=0$ at both points of $\pi^{-1}(p_0)$. Moreover, $g$ satisfies
\[
 E(g)<2E(h),\qquad g\circ\iota=-g,
\]
where $\iota$ is the nontrivial deck involution of $\pi$.
\end{lemma}

\begin{proof}
Choose a M\"obius map $M\colon\D\to\HH$ taking $p_0$ to infinity, and
write $f(x)=h(M^{-1}(x))$ for $x\in\R$. Since
$M^{-1}(z)-p_0=O(|z|^{-1})$ as $|z|\to\infty$, \eqref{eq:scalar-decay}
gives $f(x)=O(|x|^{-\delta})$ and $h(M^{-1}(\sqrt{-1}Y))=O(Y^{-\delta})$ as
$Y\to\infty$. In particular, $f\in L^1(\R)$. The function $f$ is locally
H\"older continuous because $M^{-1}$ is smooth on $\R$. By uniqueness
for the Dirichlet problem, $h\not\equiv0$ implies
$h|_{\partial\D}\not\equiv0$. Since $h(p_0)=0$ by \eqref{eq:scalar-decay}
and $M^{-1}(\R)=\partial\D\setminus\{p_0\}$, we have $f\not\equiv0$.

For $a\in\R$ and $y>0$, let
\begin{equation}\label{eq:scalar-cayley}
 A_{a,y}(z)=\frac{z-a-\sqrt{-1}y}{z-a+\sqrt{-1}y}.
\end{equation}
This map takes $\HH$ onto $\D$, sends $a+\sqrt{-1}y$ to $0$, and sends infinity
to $1$. Write $A_{a,y}(x)=e^{\sqrt{-1}\xi(x)}$ with $0<\xi(x)<2\pi$ for
$x\in\R$. With $X=x-a$ and $R=(X^2+y^2)^{1/2}$, we have
\begin{equation}\label{eq:scalar-half-angle}
 e^{\sqrt{-1}\xi(x)/2}=\frac{-X+\sqrt{-1}y}{R},\qquad \xi'(x)=\frac{2y}{R^2}.
\end{equation}
The quotient $(-X+\sqrt{-1}y)/R$ has positive imaginary part $y/R$, so
its argument lies in $(0,\pi)$. Differentiating the argument gives the formula for $\xi'(x)$
in \eqref{eq:scalar-half-angle}, and $\xi(x)$ increases from $0$ to
$2\pi$ along $\R$.

Fix $Y>0$ and set $v_Y=h\circ M^{-1}\circ A_{0,Y}^{-1}$. The map
$M^{-1}\circ A_{0,Y}^{-1}$ is a disk automorphism, so $v_Y$ is harmonic on
$\D$ and continuous on $\clD$. Since $A_{0,Y}(\sqrt{-1}Y)=0$, the
mean value property on circles of radius $r<1$, followed by $r\nearrow1$, gives
\[
 h(M^{-1}(\sqrt{-1}Y))=v_Y(0)
 =\frac{1}{2\pi}\int_0^{2\pi}v_Y(e^{\sqrt{-1}\xi})\,d\xi.
\]
Write $A_{0,Y}(x)=e^{\sqrt{-1}\xi(x)}$ for $x\in\R$. By
\eqref{eq:scalar-half-angle} with $a=0$ and $y=Y$, the angle $\xi(x)$
increases from $0$ to $2\pi$ as $x$ runs from $-\infty$ to $+\infty$, and
\[
 v_Y(e^{\sqrt{-1}\xi(x)})=h(M^{-1}(x))=f(x),
 \qquad
 \frac{d\xi}{dx}=\frac{2Y}{x^2+Y^2}.
\]
Changing variables in the boundary integral and multiplying by $\pi Y$
therefore gives
\[
 \pi Y h(M^{-1}(\sqrt{-1}Y))
 =\int_{\R}\frac{Y^2}{x^2+Y^2}f(x)\,dx.
\]
Since $\delta>1$, the left-hand side tends to zero as $Y\to\infty$.
Dominated convergence on the right yields
\begin{equation}\label{eq:scalar-mean-zero}
 \int_{\R} f(x)\,dx=0.
\end{equation}

The primitive $H(t)=\int_{-\infty}^t f(x)\,dx$ is nonconstant and tends
to zero at both ends of $\R$ by \eqref{eq:scalar-mean-zero}. We can thus
choose $a\in\R$ such that
\begin{equation}\label{eq:scalar-extremum}
 |H(a)|=\max_{t\in\R}|H(t)|>0,\qquad f(a)=H'(a)=0.
\end{equation}
For $s>0$, integration by parts with $H-H(a)$ gives
\begin{equation}\label{eq:scalar-laplace}
 \begin{aligned}
 \int_{-\infty}^{a} e^{-s(a-x)}f(x)\,dx
 &=s\int_{-\infty}^{a} e^{-s(a-x)}(H(a)-H(x))\,dx,\\
 \int_{a}^{\infty} e^{-s(t-a)}f(t)\,dt
 &=-s\int_{a}^{\infty} e^{-s(t-a)}(H(a)-H(t))\,dt.
 \end{aligned}
\end{equation}
The endpoint terms vanish because $H-H(a)$ is zero at $a$ and bounded
on $\R$. By \eqref{eq:scalar-extremum}, $H(a)-H(t)$ has the sign of
$H(a)$ wherever it is nonzero. Since $H(a)-H(t)\to H(a)\ne0$ as
$t\to\pm\infty$, the two integrals on the left in
\eqref{eq:scalar-laplace} have strictly opposite signs.

To use the Fubini theorem for the product of the two integrals in
\eqref{eq:scalar-laplace}, split the region $x<a<t$ at $t-x=\eps$.
By \eqref{eq:scalar-extremum}, $f(a)=0$. For $x<a<t$ and $t-x<\eps$,
local H\"older continuity of $f$ with exponent $0<\lambda<1$ gives, for small $\eps$,
\[
 |f(x)f(t)|\le C(a-x)^\lambda(t-a)^\lambda
 \le C(a-x)^\lambda(t-x)^\lambda.
\]
Using $f\in L^1(\R)$ on the region $t-x\ge\eps$, we obtain
\begin{align}
 \iint_{x<a<t}\frac{|f(x)f(t)|}{(t-x)^2}\,dx\,dt
 &=\iint_{\substack{x<a<t\\t-x<\eps}}
   \frac{|f(x)f(t)|}{(t-x)^2}\,dx\,dt +\iint_{\substack{x<a<t\\t-x\ge\eps}}
   \frac{|f(x)f(t)|}{(t-x)^2}\,dx\,dt \notag\\
 &\le C\int_{a-\eps}^{a}(a-x)^\lambda
         \int_a^{x+\eps}(t-x)^{\lambda-2}\,dt\,dx
       +\eps^{-2}\|f\|_{L^1(\R)}^2 \notag\\
 &\le \frac{C}{1-\lambda}\int_{a-\eps}^{a}(a-x)^{2\lambda-1}\,dx
       +\eps^{-2}\|f\|_{L^1(\R)}^2<\infty.
 \label{eq:scalar-cross-integrability}
\end{align}
Since
$\int_0^\infty s e^{-s(t-x)}\,ds=(t-x)^{-2}$ for $t>x$,
\eqref{eq:scalar-cross-integrability} justifies the Fubini theorem, and
\eqref{eq:scalar-laplace} gives the sign in
\begin{equation}\label{eq:I}
 \begin{aligned}
 I =\iint_{x<a<t}\frac{f(x)f(t)}{(t-x)^2}\,dx\,dt =\int_0^\infty s
     \left(\int_{-\infty}^{a} e^{-s(a-x)}f(x)\,dx\right)
     \left(\int_a^\infty e^{-s(t-a)}f(t)\,dt\right)\,ds<0.
 \end{aligned}
\end{equation}

With $a$ as in \eqref{eq:scalar-extremum}, take $y>0$ and set
\begin{equation}\label{eq:scalar-cover}
 \psi_y=A_{a,y}\circ M,\qquad
 \pi_y(\zeta)=\psi_y^{-1}(\zeta^2),\qquad
 f_y(\xi)=h(\psi_y^{-1}(e^{\sqrt{-1}\xi})).
\end{equation}
By \eqref{eq:scalar-cayley}, $\psi_y$ is a disk automorphism with
$\psi_y(p_0)=1$ and $\psi_y(M^{-1}(a+\sqrt{-1}y))=0$. Thus $\pi_y$ is smooth on
$\clD$ and has a single simple ramification point at $\zeta=0$. The
points over $p_0$ are $1$ and $-1$. By \eqref{eq:scalar-cover}, on the two open semicircles we have 
\[
 \pi_y(e^{\sqrt{-1}\theta})=
 \begin{cases}
 \psi_y^{-1}(e^{\sqrt{-1}(2\theta)}), & 0<\theta<\pi,\\
 \psi_y^{-1}(e^{\sqrt{-1}(2\theta-2\pi)}), & \pi<\theta<2\pi.
 \end{cases}
\]
As $\theta$ increases through either interval, the corresponding angle
runs from $0$ to $2\pi$. Since $\psi_y^{-1}$ preserves the boundary orientation and
maps $1$ to $p_0$, the restriction of $\pi_y$ to each open semicircle of $\partial \D$ is an orientation preserving diffeomorphism onto $\partial\D\setminus\{p_0\}$.
In terms of the angle $\theta$ on the domain $\pi_y$, the boundary data are defined by
\begin{equation}\label{eq:scalar-trace}
 U_y(\theta)=
 \begin{cases}
 f_y(2\theta),&0\le\theta\le\pi,\\
 -f_y(2\theta-2\pi),&\pi\le\theta\le2\pi.
 \end{cases}
\end{equation}
Since $f_y(0)=f_y(2\pi)=h(p_0)=0$ by \eqref{eq:scalar-decay}, the two
pieces agree at the endpoints and define a continuous periodic
function $U_y$ of period $2\pi$. The Poisson extension $g_y$ of $U_y$ is continuous on $\clD$.
The nontrivial deck involution of $\pi_y$ is $\iota(\zeta)=-\zeta$ by
\eqref{eq:scalar-cover}. Precomposition with $\iota$ exchanges the two
semicircles and preserves Dirichlet energy, so it suffices to use the
ordering in \eqref{eq:scalar-trace}. Since this trace satisfies
$U_y(\theta+\pi)=-U_y(\theta)$, uniqueness of the harmonic extension gives
$g_y\circ\iota=-g_y$. Moreover, since $h\circ\pi_y$ is harmonic, under the conformal change of variables for the degree-2 map $\pi_y$ we have
\begin{equation}\label{eq:scalar-lift-energy}
 E(h\circ\pi_y)=2E(h).
\end{equation}

To compare $E(g_y)$ with $E(h\circ\pi_y)$, we first express the Dirichlet
energy of a Poisson extension in terms of its boundary values. For a real
$2\pi$-periodic function $v\in L^2(0,2\pi)$, write
\[
 c_n=\frac{1}{2\pi}\int_0^{2\pi}v(\xi)e^{-\sqrt{-1}n\xi}\,d\xi,
 \qquad
 P[v](re^{\sqrt{-1}\xi})=\sum_{n\in\mathbb Z}c_n r^{|n|}e^{\sqrt{-1}n\xi}.
\]
The series can be differentiated term by term on $|z|\le\rho$ for
$0<\rho<1$. The functions $e^{\sqrt{-1}n\xi}$ are orthogonal in $L^2(0,2\pi)$, so the
cross terms between distinct modes vanish upon integration in $\xi$.
For each mode, the coefficients of $\partial_r P[v]$ and
$r^{-1}\partial_\xi P[v]$ have the same modulus, so the two terms contribute
equally to the Dirichlet energy. Thus
\begin{align*}
 E(P[v],\{|z|<\rho\}) =2\pi\sum_{n\ne0}n^2|c_n|^2\int_0^\rho r^{2|n|-1}\,dr =\pi\sum_{n\ne0}|n||c_n|^2\rho^{2|n|}.
\end{align*}
Letting $\rho\nearrow1$ gives
$E(P[v])=\pi\sum_n|n||c_n|^2$, with both sides allowed to be infinite.

To recover the same sum from the boundary values, apply Parseval's
identity to $v(\cdot+s)-v$ to obtain
\[
 \int_0^{2\pi}|v(\xi+s)-v(\xi)|^2\,d\xi
 =2\pi\sum_n|c_n|^2|e^{\sqrt{-1}ns}-1|^2.
\]
For $n\ne0$, the geometric sum gives
\[
 \frac{|e^{\sqrt{-1}ns}-1|^2}{4\sin^2(s/2)}
 =\left|\sum_{j=0}^{|n|-1}e^{\sqrt{-1}js}\right|^2.
\]
The integral of this squared sum over $0<s<2\pi$ is
$2\pi|n|$.  Divide the Parseval identity by
$4\sin^2(s/2)$ and integrate in $s$. Tonelli's theorem gives
$4\pi^2\sum_n|n||c_n|^2=4\pi E(P[v])$ for the resulting double integral.
Writing $\eta=\xi+s$ modulo $2\pi$ gives the scalar boundary formula for the Dirichlet integral~\cite[\S5, (5.1), and \S22, (22.3)]{douglas}
\begin{equation}\label{eq:scalar-fourier}
 E(P[v])=\frac{1}{4\pi}\int_0^{2\pi}\!\int_0^{2\pi}
 \frac{(v(\xi)-v(\eta))^2}{4\sin^2((\xi-\eta)/2)}\,d\xi\,d\eta.
\end{equation}

In applying \eqref{eq:scalar-fourier} to the two traces $U_y$ and $h\circ \pi_y$, we write
$(\theta,\varphi)$ for the integration variables $(\xi,\eta)$. Write
$V_y(\theta)=h(\pi_y(e^{\sqrt{-1}\theta}))$ for the trace of $h\circ\pi_y$.
Since $g_y$ and $h\circ\pi_y$ are the harmonic extensions of $U_y$ and
$V_y$, respectively, \eqref{eq:scalar-fourier} gives
\begin{align*}
 E(g_y)=\frac{1}{4\pi}\int_0^{2\pi}\!\int_0^{2\pi}
 \frac{(U_y(\theta)-U_y(\varphi))^2}{4\sin^2((\theta-\varphi)/2)}
 \,d\varphi\,d\theta,\quad
 E(h\circ\pi_y)=\frac{1}{4\pi}\int_0^{2\pi}\!\int_0^{2\pi}
 \frac{(V_y(\theta)-V_y(\varphi))^2}{4\sin^2((\theta-\varphi)/2)}
 \,d\varphi\,d\theta.
\end{align*}

By \eqref{eq:scalar-cover} and \eqref{eq:scalar-trace}, $U_y=V_y$ on the upper semicircle and $U_y=-V_y$ on the lower semicicle. The two integrands therefore agree when both angles lie in the same semicircle. For $0<\theta<\pi<\varphi<2\pi$, their difference is
\[
 \frac{(U_y(\theta)-U_y(\varphi))^2-(V_y(\theta)-V_y(\varphi))^2}
      {4\sin^2((\theta-\varphi)/2)}
 =\frac{V_y(\theta)V_y(\varphi)}{\sin^2((\theta-\varphi)/2)}.
\]
The absolute difference is symmetric in $\theta$ and $\varphi$, so its
integral on $0<\varphi<\pi<\theta<2\pi$ is the same.

On the region $0<\theta<\pi<\varphi<2\pi$, change to the base angles
$\xi=2\theta$ and $\eta=2\varphi-2\pi$. Both range over $(0,2\pi)$, and
\eqref{eq:scalar-cover} gives $V_y(\theta)=f_y(\xi)$ and
$V_y(\varphi)=f_y(\eta)$. Since $\theta=\xi/2$ and $\varphi=\eta/2+\pi$,
\begin{align*}
 4\sin^2\!\left(\frac{\theta-\varphi}{2}\right) =4\cos^2\!\left(\frac{\xi-\eta}{4}\right), \qquad d\theta\,d\varphi=\tfrac14\,d\xi\,d\eta.
\end{align*}
Including the two regions with angles in opposite semicircles, we obtain the integral
of the absolute difference of the two integrands $E(g_y)$ and $E(h\circ \pi_y)$ as
\begin{align}
 \frac{1}{2\pi}\int_0^\pi\!\int_\pi^{2\pi}
 \frac{|V_y(\theta)V_y(\varphi)|}{\sin^2((\theta-\varphi)/2)}
 \,d\varphi\,d\theta =\frac{1}{8\pi}\int_0^{2\pi}\!\int_0^{2\pi}
 |f_y(\xi)f_y(\eta)|\sec^2\bigl((\xi-\eta)/4\bigr)\,d\eta\,d\xi.
 \label{eq:scalar-absolute-difference}
\end{align}

To estimate the integral in
 \eqref{eq:scalar-absolute-difference}, use the boundary angle
of $A_{a,y}$ in \eqref{eq:scalar-cayley} to set $\xi=\xi(x)$ and
$\eta=\xi(t)$. By \eqref{eq:scalar-cover}, $f_y(\xi(x))=f(x)$ and
$f_y(\xi(t))=f(t)$. Write $T=t-a$ and $S=(T^2+y^2)^{1/2}$. Taking the
real part of $e^{\sqrt{-1}\xi(x)/2}e^{-\sqrt{-1}\xi(t)/2}$ using
\eqref{eq:scalar-half-angle}, we obtain
\[
 \cos^2\left(\frac{\xi(x)-\xi(t)}{4}\right)
 =\frac{1+\cos\bigl((\xi(x)-\xi(t))/2\bigr)}{2}
 =\frac12\left(1+\frac{XT+y^2}{RS}\right).
\]
By \eqref{eq:scalar-half-angle}, the Jacobian is
$d\xi\,d\eta=4y^2(RS)^{-2}\,dx\,dt$. Thus the change of variables in
\eqref{eq:scalar-absolute-difference} gives
\begin{equation}\label{eq:scalar-real-line-integral}
 \frac{1}{8\pi}\int_0^{2\pi}\!\int_0^{2\pi}
 |f_y(\xi)f_y(\eta)|\sec^2\!\left(\frac{\xi-\eta}{4}\right)\,d\eta\,d\xi
 =\frac{1}{\pi}\iint_{\R^2}|f(x)f(t)|k_y(x,t)\,dx\,dt,
\end{equation}
where
\begin{equation}\label{eq:scalar-kernel}
 k_y(x,t)=\frac{y^2}{RS(RS+XT+y^2)}.
\end{equation}
For $x\ne t$, multiplying the numerator and denominator in
\eqref{eq:scalar-kernel} by $RS-XT-y^2$ and recalling $R^2S^2-(XT+y^2)^2=y^2(x-t)^2$, we obtain
\begin{equation}\label{eq:scalar-kernel-limit}
 k_y(x,t)=\frac{1-(XT+y^2)/(RS)}{(x-t)^2},\qquad x\ne t.
\end{equation}
Figure~\ref{fig:scalar-kernel-limits} shows how the angle in \eqref{eq:scalar-kernel-limit} depends on the positions of $x$ and $t$ relative to $a$.

\begin{figure}[!htbp]
  \centering
  \includegraphics[width=.93\linewidth]{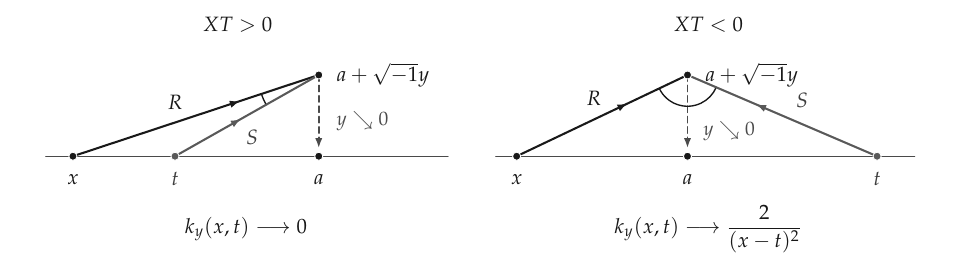}
  \caption{The arrows from $x$ and $t$ to $a+\sqrt{-1}y$ have lengths $R$ and $S$.
  The cosine of the marked angle is $(XT+y^2)/(RS)$ in
  \eqref{eq:scalar-kernel-limit}. For fixed $x,t\ne a$ with $x\ne t$,
  the angle tends to $0$ if $XT>0$ and to $\pi$ if $XT<0$ as
  $y\searrow0$.}
  \label{fig:scalar-kernel-limits}
\end{figure}

If $XT\ge0$, \eqref{eq:scalar-kernel} gives $0\le k_y\le y^2/(R^2S^2)$.
Hence
\[
 \iint_{XT\ge0}|f(x)f(t)|k_y(x,t)\,dx\,dt
 \le\left(\int_{\R}|f(x)|\frac{y}{(x-a)^2+y^2}\,dx\right)^2.
\]
For fixed $y>0$, the integral on the right is finite since
$y/((x-a)^2+y^2)\le1/y$ and $f\in L^1(\R)$. For the limit as
$y\searrow0$, split this integral at $|x-a|=\eps$ to obtain
\[
 \int_{\R}|f(x)|\frac{y}{(x-a)^2+y^2}\,dx
 \le\pi\sup_{|x-a|<\eps}|f(x)|+y\eps^{-2}\|f\|_{L^1(\R)}.
\]
Letting $y\searrow0$ and then $\eps\searrow0$ gives 0 as $f$ is continuous at $a$ and $f(a)=0$ by \eqref{eq:scalar-extremum}.

If $XT<0$, \eqref{eq:scalar-kernel-limit} gives $0\le k_y\le2/(x-t)^2$.
Also $(XT+y^2)/(RS)\to-1$ as $y\searrow0$, so $k_y\to2/(x-t)^2$.
Thus the bound in \eqref{eq:scalar-cross-integrability} is applicable on $x<a<t$ and, by symmetry, on $t<a<x$. By dominated convergence theorem and \eqref{eq:I},
each region contributes $2I$ in the limit. Hence, as $y \searrow 0$, we obtain
\[
 \iint_{XT<0} f(x)f(t)k_y(x,t)\,dx\,dt\longrightarrow 4I.
\]

The estimates on $XT\ge0$ and $XT<0$ show that the right-hand side of
\eqref{eq:scalar-real-line-integral} is finite for each $y>0$. The
boundary integral for $h\circ\pi_y$ in \eqref{eq:scalar-fourier} is finite
by \eqref{eq:scalar-lift-energy}. The boundary integral for $U_y$ is
therefore finite as well, and \eqref{eq:scalar-fourier} gives
$E(g_y)<\infty$. Together with $g_y\in C(\clD)$, this implies
$g_y\in W^{1,2}(\D)$. We can now subtract the two energy formulas to obtain
\[
 E(g_y)-2E(h)=\frac{1}{\pi}\iint_{\R^2} f(x)f(t)k_y(x,t)\,dx\,dt
 \longrightarrow\frac{4I}{\pi}<0\qquad\text{as }y\searrow0.
\]
Choose $y>0$ small enough that $E(g_y)<2E(h)$ and set $\pi=\pi_y$,
$b_0=\pi_y(0)$, and $g=g_y$, completing the proof of Lemma~\ref{lem:scalar}.
\end{proof}

\subsection{The comparison map and its area}
\label{sec:transfer}

In this subsection, we construct the comparison map $Y$ on a branched double cover of $\Omega$ with the prescribed graph boundary traces. When $u^3_-\not\equiv0$, we use isothermal coordinates to transfer the Dirichlet energy comparison to the strict area inequality \eqref{eq:drop}.

For small $r$, \eqref{eq:decay} gives
$|\nabla u^3_\pm|\le Cr^{\alpha/q}\le\eps_0$ on $\Omega$, where
$\eps_0$ is from Lemma~\ref{lem:maj}. Define
\begin{equation}\label{eq:transfer-coefficient}
 B(w)=B_{\nabla u^3_+(w),\,\nabla u^3_-(w)}.
\end{equation}
By Lemma~\ref{lem:maj}, $B$ is measurable and symmetric positive definite,
with $\det B=1$ and $\|B-I\|_{\mathrm{op}}\le Cr^{2\alpha/q}$. In particular,
$1\le\lambda_{\max}(B(w))\le1+Cr^{2\alpha/q}$. Since $\beta\theta/\pi>1$,
we can shrink $r$ so that
\begin{equation}\label{eq:dist}
 K:=\operatorname*{ess\,sup}_{w\in\Omega}\lambda_{\max}(B(w)),
 \qquad \delta=\frac{\beta\theta}{K\pi}>1.
\end{equation}

The first variation of the area of the two minimal graphs $u^3_j$ gives
\[
 \int_\Omega DW(\nabla u^3_j)\cdot\nabla\phi\,dw=0
 \qquad j=1,2,\quad  \forall \phi\in C_c^\infty(\Omega).
\]
The first identity in \eqref{eq:tangent}, with $B$ as in
\eqref{eq:transfer-coefficient}, gives
$DW(\nabla u^3_1)-DW(\nabla u^3_2)=2B\nabla u^3_-$. Subtracting the two
first variation identities gives the weak equation
$\diver(B\nabla u^3_-)=0$, namely
\begin{equation}\label{eq:transfer-weak}
 \int_\Omega\inner{B\nabla u^3_-,\nabla\phi}\,dw=0
 \qquad\text{for all }\phi\in C_c^\infty(\Omega).
\end{equation}

Let $S\colon\D\to\Omega$ be a Riemann map. Both $S$ and $S^{-1}$ are
smooth on their respective domains. The conductivity in the disk
coordinate $z=S^{-1}(w)$ is given by
\[
 ((S^{-1})_*B)(z)
 =\frac{D(S^{-1})(w)\,B(w)\,D(S^{-1})(w)^{\mathsf T}}{J_{S^{-1}}(w)},
 \qquad w=S(z),
\]
where $J_{S^{-1}}=\det D(S^{-1})$. Since $S^{-1}$ is conformal, we can
write $D(S^{-1})(w)=\rho(w)R(w)$ with $\rho(w)>0$ and $R(w)$ a rotation.
Then $J_{S^{-1}}(w)=\rho(w)^2$, so $((S^{-1})_*B)(z)=R(w)B(w)R(w)^{\mathsf T}$.
In particular, $(S^{-1})_*B$ has the same eigenvalues as $B$ at
corresponding points.

Define $\sigma=(S^{-1})_*B$ on $\D$ and $\sigma=I$ on $\C\setminus\D$.
By \eqref{eq:dist} and $\det B=1$, we have $K^{-1}I\le\sigma\le KI$ on $\mathbb{C}$ and
$\det\sigma=1$ almost everywhere. By Astala--Lassas--P\"aiv\"arinta~\cite[Lemma~3.1]{alp}, there exists a H\"{o}lder continuous quasiconformal homeomorphism
$F_0\colon\C\to\C$ with
\[
 (F_0)_*\sigma=\sqrt{\det(\sigma\circ F_0^{-1})}\,I=I.
\]
In particular, $F_0\in W^{1,2}_{\mathrm{loc}}(\C)$.
Since the image $F_0(\D)$ is a Jordan domain as $F_0$ is a homeomorphism
of the plane, we can choose a Riemann map $H\colon F_0(\D)\to\D$ normalized by
$H(F_0(0))=0$. Then $Q=H\circ(F_0|_{\D})$ maps $\D$ onto $\D$ and
satisfies $Q(0)=0$. By conformality of $H$, the product $DH\,(DH)^{\mathsf T}$ equals
$(\det DH)I$ on $F_0(\D)$. Using the previous pushforward formula gives $H_*I=I$. Since $(F_0)_*\sigma=I$ and $\sigma=(S^{-1})_*B$ on $\D$, we have
\begin{align*}
 Q_*\bigl((S^{-1})_*B\bigr) =H_*\bigl(((F_0)_*\sigma)|_{F_0(\D)}\bigr) =H_*I=I.
\end{align*}
Set $T=Q\circ S^{-1}\colon\Omega\to\D$. Since $Q=H\circ(F_0|_{\D})$
and $H$ and $S^{-1}$ are conformal, $T$ is quasiconformal. In particular,
$T\in W^{1,2}_{\mathrm{loc}}(\Omega)$. The preceding identity gives
$T_*B=I$. Evaluating the pushforward formula at $T(w)$ and multiplying by $J_T(w)$ gives
\begin{equation}\label{eq:transfer-isothermal}
 DT(w)B(w)DT(w)^{\mathsf T}=J_T(w)I
 \qquad\text{for a.e. }w\in\Omega.
\end{equation}
In particular, $(DT)^{\mathsf T}DT=J_T B^{-1}$. The eigenvalues of $B^{-1}$
are $\lambda$ and $\lambda^{-1}$, where $\lambda=\lambda_{\max}(B)$.
Hence the larger eigenvalue of $(DT)^{\mathsf T}DT$ is $J_T\lambda$, and
$\|DT\|_{\mathrm{op}}^2=J_T\lambda\le KJ_T$ almost everywhere. Thus $T$ and $Q=T\circ S$ are $K$-quasiconformal. By Carath\'eodory's theorem~\cite{caratheodory}, the Riemann maps $H$ and $S$
extend homeomorphically to the closures of their Jordan domains. Since $F_0$
is a homeomorphism of $\C$, the quasiconformal maps $Q$ and $T$ extend to
homeomorphisms $\clD\to\clD$ and $\overline\Omega\to\clD$, respectively.
Since $Q^{-1}(0)=0$, Mori's theorem~\cite[p.~157, (2)]{mori} gives
\[
 |Q^{-1}(z)-Q^{-1}(z')|\le C|z-z'|^{1/K}
 \qquad\text{for }z,z'\in\clD.
\]

To estimate $S$ at the vertex $0\in\partial\Omega$, first map $\Omega$
to $D_1^+$ by $\xi=(e^{-\sqrt{-1}\theta_0}w/r)^{\pi/\theta}$, taking
$0<\arg(e^{-\sqrt{-1}\theta_0}w)<\theta$. The map
$\eta=(1+\xi)/(1-\xi)$ has
\[
 \operatorname{Re}\eta=\frac{1-|\xi|^2}{|1-\xi|^2},\qquad
 \operatorname{Im}\eta=\frac{2\operatorname{Im}\xi}{|1-\xi|^2},
\]
so it takes $D_1^+$ onto the first quadrant of $\C$, which $\eta\mapsto\eta^2$
maps onto $\HH$. If $k\colon\D\to\HH$ is a Cayley map, we can therefore take
\[
 S(\zeta)=r e^{\sqrt{-1}\theta_0}
 \left(\frac{\sqrt{k(\zeta)}-1}{\sqrt{k(\zeta)}+1}\right)^{\theta/\pi},
\]
where $\sqrt{k(\zeta)}$ lies in the first quadrant and the outer power $\theta/\pi$ uses $0<\arg\xi<\pi$. At $\zeta_0=S^{-1}(0)$ we have $k(\zeta_0)=1$.
The fraction in the formula for $S$ defines a holomorphic extension of
$\zeta\mapsto\xi$ with a simple zero at $\zeta_0$.
In particular, $|\xi|\le C|\zeta-\zeta_0|$ near $\zeta_0$, so
\begin{equation}\label{eq:transfer-vertex}
 |S(\zeta)|\le C|\zeta-\zeta_0|^{\theta/\pi}
 \qquad\text{for }\zeta\in\clD\text{ near }\zeta_0.
\end{equation}
Since $\zeta\mapsto\xi$ is locally Lipschitz and the power
$\xi\mapsto\xi^{\theta/\pi}$ is $\theta/\pi$-H\"older on
$\overline{D_1^+}$, $S$ is $\theta/\pi$-H\"older on a neighborhood of
$\zeta_0$ in $\clD$. At the preimage of $\xi=-1$, $k$ has a simple zero and the fraction
$(\sqrt{k}-1)/(\sqrt{k}+1)$ is holomorphic as a function of $\sqrt{k}$
near $\sqrt{k}=0$. At the preimage of $\xi=1$, $1/k$ has a simple zero
and we write the fraction as
\[
 \frac{1-1/\sqrt{k}}{1+1/\sqrt{k}}.
\]
The power $\xi\mapsto\xi^{\theta/\pi}$ extends holomorphically near
$\xi=\pm1$. The $1/2$-H\"older bounds for $\sqrt{k}$ and $1/\sqrt{k}$
at the respective preimages therefore give $1/2$-H\"older bounds for
$S$ on neighborhoods of these two points in $\clD$. Away from these
three corner preimages, $S$ is locally Lipschitz on $\clD$ as it
extends holomorphically across the remaining boundary arcs. Taking a finite covering of $\clD$, we conclude that $S\in C^{0,\min\{\theta/\pi,1/2\}}(\clD)$.

Set $p_0=T(0)=Q(\zeta_0)\in\partial\D$ and $h=u^3_-\circ T^{-1}$.
Since $u^3_-$ is Lipschitz and $T^{-1}=S\circ Q^{-1}$ is H\"older
continuous, $h$ is H\"older continuous on $\clD$. By \eqref{eq:decay},
\eqref{eq:transfer-vertex} and the bound for $Q^{-1}$, for $\zeta$ near $p_0$
we have
\begin{align*}
 |h(\zeta)| \le C|S(Q^{-1}(\zeta))|^\beta
 \le C|Q^{-1}(\zeta)-\zeta_0|^{\beta\theta/\pi}  \le C|\zeta-p_0|^\delta,
\end{align*}
with $\delta>1$ as in \eqref{eq:dist}.

Since $u^3_-$ is Lipschitz and $T^{-1}$ is quasiconformal,
$h\in W^{1,2}_{\mathrm{loc}}(\D)$ and we have 
\[
 \nabla u^3_-(w)=DT(w)^{\mathsf T}\nabla h(T(w))
 \qquad\text{for a.e. }w\in\Omega.
\]
Using \eqref{eq:transfer-isothermal} and the quasiconformal change of
variables given by Ahlfors--Bers~\cite[Theorem~5(vi)]{ahlforsbers}, we obtain
\begin{equation}\label{eq:transfer-h-energy}
 \begin{aligned}
 2E(h) =\int_\Omega |\nabla h(T(w))|^2 J_T(w)\,dw =\int_\Omega\inner{B\nabla u^3_-,\nabla u^3_-}\,dw.
 \end{aligned}
\end{equation}
The last integral is finite by \eqref{eq:decay} and the bound on $B$. Thus $h\in W^{1,2}(\D)$ as $h$ is bounded on $\clD$. Since $T\in W^{1,2}_{\mathrm{loc}}(\Omega)$, for each $\phi\in C_c^\infty(\D)$ we have
$\phi\circ T\in W^{1,2}_{\mathrm{loc}}(\Omega)$, with
$\nabla(\phi\circ T)(w)=DT(w)^{\mathsf T}\nabla\phi(T(w))$ almost everywhere.
Since $\phi\circ T$ is supported in the compact set
$T^{-1}(\operatorname{supp}\phi)\subset\Omega$, it belongs to
$W^{1,2}(\Omega)$. We may plug it in \eqref{eq:transfer-weak} by smooth approximation. By
\eqref{eq:transfer-isothermal}, we have
\[
 \int_\D \nabla h\cdot\nabla\phi\,d\zeta
 =\int_\Omega\inner{B\nabla u^3_-,\nabla(\phi\circ T)}\,dw=0.
\]
Thus, by Weyl's lemma and continuity of $h$, $h$ is a smooth harmonic function on $\D$.

Suppose $u^3_-\not\equiv0$, so $h\not\equiv0$. Lemma~\ref{lem:scalar}
gives a holomorphic double cover $\Pi\colon\D\to\D$ with branch value
$b_0\in\D$ and a harmonic function $g\in C(\clD)\cap W^{1,2}(\D)$ with
$E(g)<2E(h)$. The trace of $g$ equals $h\circ\Pi$ on one component of
$\partial\D\setminus\Pi^{-1}(p_0)$ and $-h\circ\Pi$ on the other.
Moreover, $g=0$ on $\Pi^{-1}(p_0)$. If $\iota$ is the nontrivial deck
involution of $\Pi$, then $g\circ\iota=-g$ by Lemma~\ref{lem:scalar}.

\begin{figure}[!htbp]
  \centering
  \includegraphics[width=.91\linewidth]{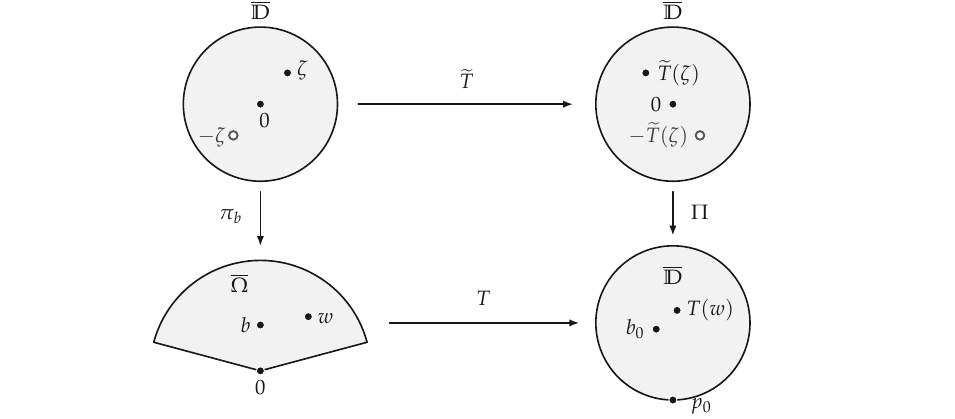}
  \caption{The lift $\widetilde T$ in \eqref{eq:transfer-lift}. For
  $w=\pi_b(\zeta)\ne b$, the filled and open points labelled $\zeta$
  and $-\zeta$ map to the corresponding points over $T(w)$.
  The upper disks have ramification point $0$. On the base, $T(b)=b_0$
  and $T(0)=p_0$.}
  \label{fig:transfer-lift}
\end{figure}

Set $b=T^{-1}(b_0)$. Let $\pi_b\colon\D\to\Omega$ be the holomorphic
double cover given by $\pi_b(\zeta)=c_b^{-1}(\zeta^2)$, where
$c_b\colon\Omega\to\D$ is conformal and $c_b(b)=0$. Choose $\Pi$ as
in \eqref{eq:scalar-cover} and write $c_{b_0}=\psi_y$, so that
$\Pi(\zeta)=c_{b_0}^{-1}(\zeta^2)$ and $c_{b_0}(b_0)=0$. Then
$F=c_{b_0}\circ T\circ c_b^{-1}$ is an orientation preserving homeomorphism of the closed disk fixing $0$. On $\D$, $F$ is $K$-quasiconformal by conformality of $c_{b_0}$ and
$c_b^{-1}$, so $F\in W^{1,2}_{\mathrm{loc}}(\D)$.
For $z\ne0$, $F(z)/z$ has winding number zero along each positively
oriented circle centered at $0$, since $F(z)$ and $z$ each wind once.
Choose a continuous square root of $F(z)/z$ on $\clD\setminus\{0\}$ and define a lift of $T$ by
\[
 \widetilde T(\zeta)=\zeta\sqrt{\frac{F(\zeta^2)}{\zeta^2}}
 \quad(\zeta\ne0),\qquad \widetilde T(0)=0.
\]
The equality $|\widetilde T(\zeta)|^2=|F(\zeta^2)|$ gives continuity of $\widetilde T$ at
$0$. For $z\in\clD\setminus\{0\}$, let $\zeta$ be a square root of $z$. The
formula for $\widetilde T$ gives $\widetilde T(\zeta)^2=F(z)$ and
$\widetilde T(-\zeta)=-\widetilde T(\zeta)$. Since $F(z)\ne0$, the two
images $\widetilde T(\zeta)$ and $\widetilde T(-\zeta)$ are precisely the
two distinct square roots of $F(z)$. Since $F$ is bijective and
$\widetilde T(0)=0$, the continuous map $\widetilde T$ is a bijection of $\clD$ onto itself
and hence a homeomorphism. It satisfies
\begin{equation}\label{eq:transfer-lift}
 \Pi\circ\widetilde T=T\circ\pi_b,
 \qquad \widetilde T\circ\iota_b=\iota\circ\widetilde T,
\end{equation}
where $\iota_b(\zeta)=-\zeta$ is the nontrivial deck involution of $\pi_b$.
Figure~\ref{fig:transfer-lift} shows how $\widetilde T$ maps the two points over $w$ to the two points over $T(w)$.

For each $\zeta\in\D\setminus\{0\}$, the maps $\pi_b$ and $\Pi$ are locally conformal
at $\zeta$ and $\widetilde T(\zeta)$, respectively. By the first identity
in \eqref{eq:transfer-lift}, $\widetilde T$ is $K$-quasiconformal on $\D\setminus\{0\}$. At almost every point of $\D\setminus\{0\}$, write $s_1^2$ and $s_2^2$
for the eigenvalues of $(D\widetilde T)^{\mathsf T}D\widetilde T$,
with $s_1\ge s_2>0$. Then $J_{\widetilde T}=s_1s_2$, and
$K$-quasiconformality gives $s_1/s_2\le K$. Consequently
\[
 |D\widetilde T|^2=s_1^2+s_2^2
 =\left(\frac{s_1}{s_2}+\frac{s_2}{s_1}\right)J_{\widetilde T}
 \le (K+K^{-1})J_{\widetilde T}.
\]
Integrating and changing variables by the quasiconformal homeomorphism
$\widetilde T\colon\D\setminus\{0\}\to\D\setminus\{0\}$ gives
\begin{align*}
 \int_{\D\setminus\{0\}}|D\widetilde T|^2\,d\zeta \le (K+K^{-1})\int_{\D\setminus\{0\}}J_{\widetilde T}\,d\zeta =(K+K^{-1})|\widetilde T(\D\setminus\{0\})|=(K+K^{-1})\pi.
\end{align*}

To extend the weak derivative of $\widetilde T$ across $0$, let $v$ be a bounded map on
$\D$ with $v\in W^{1,2}_{\mathrm{loc}}(\D\setminus\{0\})$ and
$Dv\in L^2(\D\setminus\{0\})$. For $0<\eps<1$, set
\begin{equation}\label{eq:transfer-log-cutoff}
 \eta_\eps(\zeta)=
 \begin{cases}
  0 & |\zeta|\le\eps^2,\\
  \dfrac{\log(|\zeta|/\eps^2)}{\log(1/\eps)} & \eps^2<|\zeta|<\eps,\\
  1 & |\zeta|\ge\eps.
 \end{cases}
\end{equation}
On $\eps^2<|\zeta|<\eps$, we have
$|\nabla\eta_\eps(\zeta)|=(|\zeta|\log(1/\eps))^{-1}$ and the gradient
vanishes outside this annulus, which implies 
\[
 \int_\D |\nabla\eta_\eps|^2\,d\zeta
 =\frac{2\pi}{\log^2(1/\eps)}\int_{\eps^2}^{\eps}\frac{d\rho}{\rho}
 =\frac{2\pi}{\log(1/\eps)}\longrightarrow0.
\]
For $\phi\in C_c^\infty(\D)$, the product $\eta_\eps\phi$ is a valid
weak test function on the punctured disk by smooth approximation.
For each coordinate derivative $\partial_j$, we obtain
\[
 \int_\D v\eta_\eps\,\partial_j\phi\,d\zeta
 =-\int_\D (\partial_jv)\eta_\eps\phi\,d\zeta
  -\int_\D v\phi\,\partial_j\eta_\eps\,d\zeta.
\]
The last term is bounded by
$\|v\|_{L^\infty}\|\phi\|_{L^2}\|\nabla\eta_\eps\|_{L^2}$ and tends to
zero. Since $0\le\eta_\eps\le1$ and $\eta_\eps\to1$ away from $0$, dominated
convergence in the other terms shows that $Dv$ is also the weak
derivative on $\D$. Thus $v\in W^{1,2}(\D)$.
Applying this argument to $\widetilde T$ gives
$\widetilde T\in W^{1,2}(\D)$. Its distortion bound holds almost
everywhere on $\D$, so $\widetilde T$ is $K$-quasiconformal on $\D$.

Define the comparison map by
\begin{equation}\label{eq:newmap}
 Y=(\pi_b,\ u^3_+\circ\pi_b+g\circ\widetilde T),
\end{equation}
identifying the horizontal plane with $\C$. Since $\Omega$ is a Jordan
domain, Carath\'eodory's theorem extends $c_b^{-1}$ continuously to
$\clD$. Hence $\pi_b(\zeta)=c_b^{-1}(\zeta^2)$ extends continuously to a map
$\clD\to\overline\Omega$. Since the height $u^3_+$ is continuous on
$\overline\Omega$, and $g$ and $\widetilde T$ are continuous on $\clD$, both components in \eqref{eq:newmap} are therefore continuous on $\clD$.

On each component of
$\partial\D\setminus\pi_b^{-1}(0)$, the trace of $g$ and
\eqref{eq:transfer-lift} give
\[
 g\circ\widetilde T=\pm h\circ\Pi\circ\widetilde T
 =\pm u^3_-\circ\pi_b,
\]
with opposite signs on the two boundary components. Each component projects by
$\pi_b$ onto $\partial\Omega\setminus\{0\}$ with positive orientation.
Since $u^3_++u^3_-=u^3_1$ and $u^3_+-u^3_-=u^3_2$, the boundary of $Y$
traverses each graph boundary once, with both traversals starting and
ending at $u(p)=0$. Replacing $\widetilde T$ by $\iota\circ\widetilde T$
interchanges the two traversals.

Since $s_1^2$ is the larger eigenvalue of
$(D\widetilde T)^{\mathsf T}D\widetilde T$, we have almost everywhere
\[
 \|D\widetilde T\|_{\mathrm{op}}^2=s_1^2
 =\frac{s_1}{s_2}J_{\widetilde T}\le KJ_{\widetilde T},
\]
which implies
\[
 |\nabla(g\circ\widetilde T)|^2
 =|(D\widetilde T)^{\mathsf T}(\nabla g\circ\widetilde T)|^2
 \le KJ_{\widetilde T}(|\nabla g|^2\circ\widetilde T).
\]
Integrating and changing variables by $\widetilde T\colon\D\to\D$ gives
\[
 E(g\circ\widetilde T)
 \le \frac K2\int_\D |\nabla g(\widetilde T(\zeta))|^2
                     J_{\widetilde T}(\zeta)\,d\zeta
 =K E(g).
\]

Since $\pi_b$ is a holomorphic double cover, we obtain
\begin{align*}
 E(\pi_b)
 =2|\Omega|,\qquad
 E(u^3_+\circ\pi_b)  =2E(u^3_+,\Omega).
\end{align*}
By \eqref{eq:decay}, $\nabla u^3_+$ is bounded, so all three energies
are finite. Since $Y$ is bounded on $\clD$, we conclude that
$Y\in W^{1,2}(\D,\R^3)$.

Let $U\subset\Omega\setminus\{b\}$ be an open set doubly covered by
$\pi_b$. Choose an inverse branch $\sigma\colon U\to\D$ and set
$v=(g\circ\widetilde T)\circ\sigma$ on $U$. The other branch is
$\iota_b\circ\sigma$, so \eqref{eq:transfer-lift} and $g\circ\iota=-g$
give $-v$ on that sheet. By the first identity in \eqref{eq:transfer-lift},
the two corresponding inverse branches of $\Pi$ on $T(U)$ pull $g$ back
to $v\circ T^{-1}$ and $-v\circ T^{-1}$. Changing variables by the two
conformal branches and then using \eqref{eq:transfer-isothermal} gives
\begin{align*}
 E(g,\Pi^{-1}(T(U)))
 =\int_{T(U)}|\nabla(v\circ T^{-1})|^2\,d\zeta =\int_U\inner{B\nabla v,\nabla v}\,dw.
\end{align*}
We sum the identity over a measurable partition of
$\Omega\setminus\{b\}$ by subsets of doubly covered open sets to obtain
\begin{equation}\label{eq:energytransport}
 E(g)=\int_\Omega\inner{B\nabla v,\nabla v}\,dw.
\end{equation}

Let $A(u^3_j)=\int_\Omega W(\nabla u^3_j)\,dw$ be the area of the $j$th
graph. By \eqref{eq:newmap}, the two sheets of $Y$ over $U$ have heights
$u^3_++v$ and $u^3_+-v$. Changing variables by the inverse branches of
$\pi_b$ gives the sum of the two graph areas in the $w$-coordinate.
Applying the second inequality in \eqref{eq:tangent} with $B$ from
\eqref{eq:transfer-coefficient}, and then using
\eqref{eq:transfer-h-energy} and \eqref{eq:energytransport}, we obtain
\begin{equation}\label{eq:drop}
 \begin{aligned}
 A(Y)-A(u^3_1)-A(u^3_2)
 &=\int_\Omega\Bigl(F_{\nabla u^3_+}(\nabla v)
                 -F_{\nabla u^3_+}(\nabla u^3_-)\Bigr)\,dw\\
 &\le \int_\Omega\left(\inner{B\nabla v,\nabla v}
                  -\inner{B\nabla u^3_-,\nabla u^3_-}\right)\,dw\\
 &=E(g)-2E(h)<0.
 \end{aligned}
\end{equation}

If $u^3_-\equiv0$, then $u^3_1=u^3_2=u^3_+$. Choose any $b\in\Omega$
and set $Y=(\pi_b,u^3_+\circ\pi_b)$. On each component of
$\partial\D\setminus\pi_b^{-1}(0)$, the projection $\pi_b$ traverses
$\partial\Omega\setminus\{0\}$ once in the positive direction. Thus
$Y|_{\partial\D}$ traverses the boundary of the graph of $u^3_+$ twice.
Each $w\in\Omega\setminus\{b\}$ has two preimages under $\pi_b$, and
both local sheets have area density $W(\nabla u^3_+)$ in the coordinate of $w$.
Consequently, we obtain
\[
 A(Y)=2\int_\Omega W(\nabla u^3_+)\,dw=A(u^3_1)+A(u^3_2).
\]

\vspace{0.5cm}
\section{Construction of a comparison disk}
\label{sec:surgery}

In this section, we cut and reglue the two sectors $S_1$ and $S_2$ of
the original parameter disk so that they form the branched double cover
$\pi_b\colon\D\to\Omega$ from Subsection~\ref{sec:transfer}. On this cover
we use the comparison map $Y$, and on the rest of the disk we keep $u$.
If $u^3_1\not\equiv u^3_2$ on $\Omega$, the replacement strictly decreases
the area by \eqref{eq:drop} and \eqref{eq:sewn}, contradicting the area
minimality of $u$.
When the graph heights $u^3_1$ and $u^3_2$ coincide on $\Omega$, the
replacement preserves area. We then reparametrize the comparison map
to make it weakly conformal and complete the proof of
Theorem~\ref{thm:main}.

\subsection{The topology and the boundary traversals}\label{sec:surgery-topology}

In this subsection, we cut the parameter disk along a slit in each of
$S_1$ and $S_2$ and reglue the banks to obtain a disk $M$. We identify
the union of the reglued sectors with the domain of the branched double
cover $\pi_b$ from Subsection~\ref{sec:transfer} and match the boundary
traversals of $Y$ to $\partial S_1$ and $\partial S_2$.

\begin{figure}[!htbp]
  \centering
  \includegraphics[width=.95\linewidth]{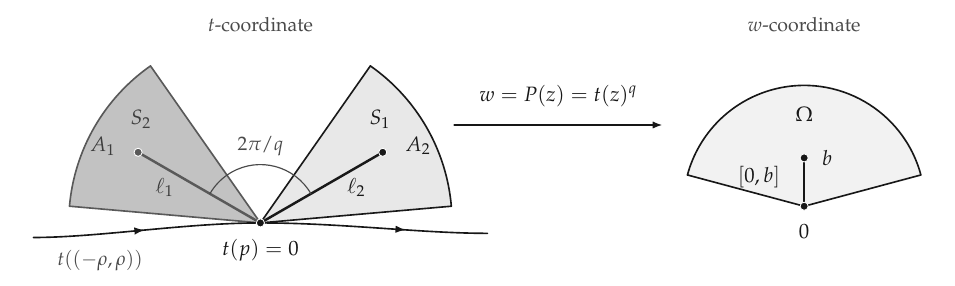}
  \caption{The two lifts of $[0,b]$ under $P=t^q$, drawn for $q=3$.
  The rays containing $\ell_2\subset S_1$ and $\ell_1\subset S_2$ differ in
  argument by $2\pi/q$. The arrows on the image of the real diameter
  under $t$ give the positive boundary direction at $p$.}
  \label{fig:surgery-slits}
\end{figure}

Let $b\in\Omega$ be the branch value of $\pi_b$. We use the inverse
branches $z_j$ of $P=t^q$ from Section~\ref{sec:geometry} to lift
$[0,b]\subset\overline\Omega$. In the $z$-coordinate, set
\[
 \ell_1=z_2([0,b]),\quad \ell_2=z_1([0,b]),\qquad
 A_1=z_2(b),\quad A_2=z_1(b).
\]
The slits start at $p$ and end at the distinct points $A_1,A_2\in\D$.
For $0\le\lambda\le1$, we have
$t(z_j(\lambda b))=\lambda^{1/q}t(z_j(b))$. Thus $\ell_1\subset S_2$
and $\ell_2\subset S_1$ are radial in the $t$-coordinate, with the
argument of $\ell_1$ larger by $2\pi/q$.
Figure~\ref{fig:surgery-slits} shows the two lifts in the $t$-coordinate.

Cutting along $\ell_1$ and $\ell_2$ produces three copies of $p$. We denote by $P_0$ and $P_2$ the limits in the cut disk along the
original boundary as $z\to0^-$ and $z\to0^+$, respectively. The third
copy $P_1$ is the vertex of the sector between $\ell_1$ and $\ell_2$. Write $\ell_i^+$ for the bank traversed toward $A_i$ and
$\ell_i^-$ for the bank traversed away from $A_i$. The positive boundary
traversal follows the arc inherited from $\partial\D\setminus\{p\}$
from $P_2$ to $P_0$. Since $\ell_1$ has the larger $t$-argument, the
traversal continues from $P_0$ along $\ell_1^+$. The order along the banks is
\begin{equation}\label{eq:surgery-bank-order}
 P_0\xrightarrow{\ell_1^+} A_1\xrightarrow{\ell_1^-}P_1
    \xrightarrow{\ell_2^+} A_2\xrightarrow{\ell_2^-}P_2.
\end{equation}
At $P_0,P_1,P_2$ the cut surface has sector neighborhoods. At an
interior point of a slit bank it has a half disk neighborhood, and at
$A_i$ a slit neighborhood of angle $2\pi$ opens to a half disk. The cut
surface is connected because a path crossing $\ell_i$ can detour around
$A_i$. In a cell decomposition containing $\ell_1$ and $\ell_2$ as
edges, cutting adds two edges and two vertices and leaves the faces
unchanged. The cut surface therefore has Euler characteristic one and
is a disk.

\begin{figure}[!htbp]
  \centering
  \includegraphics[width=.95\linewidth]{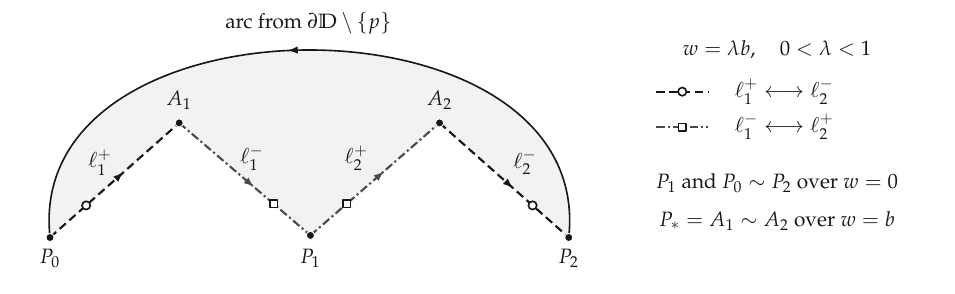}
  \caption{The boundary traversal in \eqref{eq:surgery-bank-order}.
  For a fixed $w=\lambda b$ with $0<\lambda<1$, the two open circles
  mark points identified by the first pairing in \eqref{eq:surgery-bank-pairing},
  and the two open squares mark points identified by the second.}
  \label{fig:surgery-banks}
\end{figure}

On each closed bank $\ell_i^\pm$, the coordinate $w=P(z)=t(z)^q$
gives a homeomorphism onto $[0,b]$, where $z$ denotes the coordinate of
the corresponding point before the cut. For each $w=\lambda b$ with
$0\le\lambda\le1$, we identify the points over $w$ on each of the
following pairs of banks:
\begin{equation}\label{eq:surgery-bank-pairing}
 \ell_1^+\longleftrightarrow\ell_2^-,\qquad
 \ell_1^-\longleftrightarrow\ell_2^+.
\end{equation}
In the traversal \eqref{eq:surgery-bank-order}, $\lambda$ increases on
each $\ell_i^+$ and decreases on each $\ell_i^-$. Thus both pairing
homeomorphisms reverse the boundary orientation. At $w=0$, the first
pairing identifies $P_0$ with $P_2$, and the second pairs $P_1$ with
itself. At $w=b$, both pairings identify $A_1$ with $A_2$.
Figure~\ref{fig:surgery-banks} shows the boundary traversal and the two bank pairings.

The interiors of the banks are disjoint, and the pairings agree at
their shared endpoints. Therefore the equivalence relation consists of
the diagonal and the graphs of the two pairing homeomorphisms and their
inverses. The diagonal is closed, and each graph is compact because
its domain is a closed bank. The relation is thus closed in the product
of the cut disk with itself. Thus, the quotient $M$ is compact Hausdorff and
is connected as a continuous image of the cut disk.

Away from the endpoints, each bank identification in
\eqref{eq:surgery-bank-pairing} joins two half disk neighborhoods to
form a disk. At $P_*=A_1\sim A_2$, the two slit neighborhoods are joined
along both pairs of banks, so $P_*$ also has a disk neighborhood. At
$P_1$, identifying the two sides $\ell_1^-$ and $\ell_2^+$ of its sector
similarly gives a disk neighborhood. At $P_0\sim P_2$, identifying
$\ell_1^+$ and $\ell_2^-$ joins the two sectors there into a half disk.
The remaining sides belong to the arc inherited from
$\partial\D\setminus\{p\}$. Thus $M$ is an oriented surface with one
boundary component, obtained by joining the endpoints of that arc.
The identifications remove two edges and merge two pairs of distinct
vertices, leaving the Euler characteristic equal to one. Hence $M$ is
a topological disk.

Under $w=P(z)$, cutting $S_1$ along $\ell_2$ and $S_2$ along $\ell_1$
gives two copies of $\overline\Omega$ cut along $[0,b]$. The pairings
\eqref{eq:surgery-bank-pairing} join opposite banks of $[0,b]$ in the
two copies. Since paired points have the same $w$-value, $P$ descends to
a projection on the union of the reglued sectors. This projection is
locally one-to-one away from $P_*$. Near $P_*$, the two slit neighborhoods give total angle $4\pi$
around $w=b$, and the local coordinate $\eta$ satisfies $w=b+\eta^2$ and $\eta(P_*) =0$. The two points over $w=0$ are $P_1$ and $P_0\sim P_2$.
Fix $0<\varepsilon<|b|$. In the cut sector containing $\ell_i$, the
pieces with $|w|\le\varepsilon$ bounded by $\ell_i^+$ and $\ell_i^-$
project homeomorphically onto the parts of
$\overline\Omega\cap\{|w|\le\varepsilon\}$ with $\arg w\ge\arg b$ and
$\arg w\le\arg b$, respectively, with $0$ included in both parts. By \eqref{eq:surgery-bank-pairing}, we join the pieces bounded by
$\ell_1^-$ and $\ell_2^+$ at $P_1$, and those bounded by $\ell_1^+$ and
$\ell_2^-$ at $P_0\sim P_2$. In each case the two projected parts meet
along $[0,b]\cap\{|w|\le\varepsilon\}$, and the pairing identifies the
two preimages of each point of this segment. Thus the projection maps
each pair of glued pieces homeomorphically onto
$\overline\Omega\cap\{|w|\le\varepsilon\}$ and is a local homeomorphism
at both points over $0$.

To identify the projection with $\pi_b(\zeta)=c_b^{-1}(\zeta^2)$ from
Subsection~\ref{sec:transfer}, consider $c_b(w)/(w-b)$, with value
$c_b'(b)$ at $b$. This quotient is holomorphic on $\Omega$ and extends
continuously and without zeros to $\overline\Omega$, since $c_b$ extends
to a homeomorphism of the closures and has a simple zero at $b$.
Since $\overline\Omega$ is convex, the quotient has a continuous square
root $k_b$ on $\overline\Omega$ that is holomorphic on $\Omega$. Choose the branch of $\sqrt{w-b}$ cut along the ray from $b$ through $0$. On the copy from $S_1$, set $\zeta=k_b(w)\sqrt{w-b}$, and on the
copy from $S_2$, set $\zeta=-k_b(w)\sqrt{w-b}$. The limits of
$\sqrt{w-b}$ on the two sides of the cut differ by a sign. Thus the
assigned $\zeta$-values agree on the banks paired in
\eqref{eq:surgery-bank-pairing}, including at $w=0$. They also tend to
$0$ at $P_*$. These assignments give a continuous map from the union of
the reglued sectors to $\clD$ with $\zeta^2=c_b(w)$. Over every
$w\ne b$, the two points are sent to the two distinct square roots of
$c_b(w)$, and $P_*$ is sent to $0$. The map is therefore bijective and hence a homeomorphism. In the $\zeta$-coordinate the
projection is $w=\pi_b(\zeta)$. In the interior away from $P_*$, the identification is conformal in the
$w$-charts because each local branch of $\sqrt{c_b(w)}$ has nonzero
derivative. At $P_*$, where $w=b+\eta^2$, we have
$\zeta=\pm\eta\,k_b(b+\eta^2)$ with a fixed sign, by continuity on the
punctured neighborhood. Since $k_b(b)\ne0$, the change from $\eta$ to
$\zeta$ is also conformal at $P_*$. We call the two arcs in $M$ obtained from the bank pairs in
\eqref{eq:surgery-bank-pairing} the \emph{seams}. The common image of
$\ell_1^+$ and $\ell_2^-$ runs from $P_0\sim P_2$ to $P_*$, and the
common image of $\ell_1^-$ and $\ell_2^+$ runs from $P_1$ to $P_*$.
In the $\zeta$-coordinate, the interiors of these arcs are the two
connected components of $\pi_b^{-1}(\{\lambda b:0<\lambda<1\})$.

\begin{figure}[!htbp]
  \centering
  \includegraphics[width=.93\linewidth]{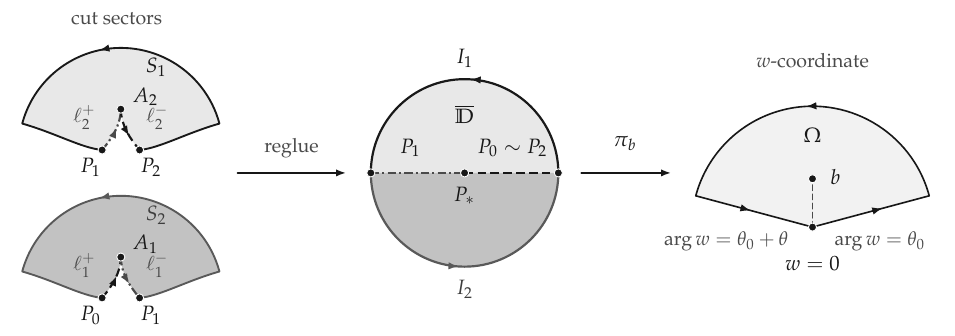}
  \caption{The two cut sectors, the covering disk, and the base sector.
  The two labels $P_1$ denote the same point of the cut disk.
  In the middle disk, the two arcs drawn from $P_*$ to the boundary
  together form $\pi_b^{-1}([0,b])$. Each $I_j$ traverses
  $\partial\Omega\setminus\{0\}$ once under $\pi_b$ in the positive direction.}
  \label{fig:surgery-traversals}
\end{figure}

Under the identification with the closed covering disk $\clD$, let
$I_j\subset\partial\D$ be the image of $\partial S_j\setminus\{p\}$
after cutting and regluing. The sets $I_1$ and $I_2$ are the two open
arcs of $\partial\D\setminus\pi_b^{-1}(0)$. The restriction
$\pi_b|_{I_j}$ maps $I_j$ homeomorphically onto
$\partial\Omega\setminus\{0\}$ and preserves the boundary orientation.
By \eqref{eq:surgery-bank-order} and \eqref{eq:surgery-bank-pairing},
$I_1$ runs from $P_0\sim P_2$ to $P_1$ in the positive direction,
and $I_2$ runs from $P_1$ to $P_0\sim P_2$. Under $\pi_b$, each
traversal starts along the side $\arg w=\theta_0$ and returns to
$w=0$ along the side $\arg w=\theta_0+\theta$.
Figure~\ref{fig:surgery-traversals} follows $I_1$ and $I_2$ from the cut sectors to $\partial\Omega$.

If $u^3_-\not\equiv0$, replacing $\widetilde T$ by
$\iota\circ\widetilde T$ in \eqref{eq:newmap} changes the sign of
$g\circ\widetilde T$ and exchanges the two graph traces. We choose the
lift so that the graph of $u^3_1$ is traced on $I_1$. The boundary
values are then
\[
 Y(\zeta)=(\pi_b(\zeta),u^3_j(\pi_b(\zeta)))
 \qquad \text{for }\zeta\in I_j,\quad j=1,2.
\]
If $u^3_-\equiv0$, the same formula holds for $Y$ because
$u^3_1=u^3_2=u^3_+$.

\subsection{Sobolev regularity of the replacement}\label{sec:surgery-sobolev}

In this subsection, we construct smooth coordinates on $M$.
On either side of each seam formed by \eqref{eq:surgery-bank-pairing},
the change from the original parameter coordinates to the new
coordinates and its inverse have uniformly bounded derivatives, even
arbitrarily close to the slit endpoints. We use these bounds to
control the Dirichlet energy when joining $Y$ to $u$ along the images
of $\partial S_1$ and $\partial S_2$, and show that the resulting map
belongs to $\A(\Gamma)$.

Away from the slit banks and their endpoints, use the original
parameter charts of $M$. At an interior point of either pair in
\eqref{eq:surgery-bank-pairing}, equal values of $w=t^q$ identify the
$t$-coordinates by a rotation through $2\pi/q$. In the original $z$-coordinate, define the chart by
$z\mapsto e^{-2\pi\sqrt{-1}/q}t(z)$ on the half neighborhood along the
$\ell_1$ bank and by $z\mapsto t(z)$ on the half-neighborhood along the
paired $\ell_2$ bank. These maps agree at paired points, and the opposite
boundary orientations in \eqref{eq:surgery-bank-pairing} place the interiors
of the two half neighborhoods on opposite sides of the common ray.
Since the rotation is an isometry and $\|Dt-I\|_{\mathrm{op}}\le 1/2$ by
Section~\ref{sec:geometry}, the derivative of each coordinate map has
operator norm at most $3/2$, and the derivative of its inverse has
operator norm at most $2$.

Choose $\eps>0$ so that $\overline{B_\eps(b)}\subset\Omega$. A neighborhood of
$P_*=A_1\sim A_2$ in $M$ consists of two copies of $B_\eps(b)$, each cut
along $[0,b]\cap B_\eps(b)$, with opposite banks identified by
\eqref{eq:surgery-bank-pairing}. One positive turn about $b$ in the
$w$-plane passes from one copy of $B_\varepsilon(b)$ to the other, and a second turn returns
to the starting point. Thus the polar angle of $w-b$ on the joined disks
is defined modulo $4\pi$. Using the coordinate $\eta$ with $w-b=\eta^2$
and $\eta(P_*)=0$, write $r=|w-b|$ and $\vartheta=2\arg\eta$ modulo
$4\pi$. We take $|\eta|\eta$ as the chart at $P_*$. Its polar form is
\[
 |\eta|\eta=r e^{\sqrt{-1}\vartheta/2}.
\]
This sends each circle of radius $r$ around the joined tips once around
the circle of radius $r$ in the chart and sends $P_*$ to $0$. For $i=1,2$, the projection $w=t^q$ identifies the slit neighborhood of
$A_i$ with one copy of $B_\eps(b)$ cut along $[0,b]\cap B_\eps(b)$.
On each copy the change from $t$ to $w-b$ is $t\mapsto t^q-b$, whose
derivative has modulus $q|t|^{q-1}=q|w|^{(q-1)/q}$. Taking
$\eps<|b|/2$ gives $|b|/2<|w|<3|b|/2$ throughout $B_\eps(b)$.
Thus the coordinate change and its inverse branch have uniformly bounded
derivatives on each of the two slit disks.

At $P_1$, put $\varphi=(\arg b)/q$, with $\theta_0<\arg b<\theta_0+\theta$.
In the $t$-plane the sector of the cut disk at $P_1$ is bounded by
$\ell_2^+$ on the ray $\arg t=\varphi$ and $\ell_1^-$ on the ray
$\arg t=\varphi+2\pi/q$. By \eqref{eq:surgery-bank-pairing}, the two
banks are identified at equal values of $w=t^q$, hence at equal radii
$|t|$. On this sector use the coordinate
\[
 t\longmapsto |t|e^{\sqrt{-1}q(\arg t-\varphi)},
 \qquad P_1\longmapsto 0.
\]
For each fixed $|t|>0$ sufficiently small, the new angle increases from
$0$ on $\ell_2^+$ to $2\pi$ on $\ell_1^-$. The two angle endpoints therefore have the same image, and the intervening arc covers the rest of the
circle exactly once. Since the radius $|t|$ is preserved, the map gives
a homeomorphism from the glued neighborhood of $P_1$ onto a disk.

To construct the chart at $P_0\sim P_2$, recall from
\eqref{eq:surgery-bank-order} that the sector at $P_2$ lies between the
positive real $z$-axis and $\ell_2^-$, and the sector at $P_0$ lies between
$\ell_1^+$ and the negative real $z$-axis. To make all four sides straight
in one coordinate, choose a smooth $\chi\colon[0,\pi]\to[0,1]$ that vanishes
near $0$ and $\pi$ and equals one on a neighborhood of the two closed
angular intervals of $t(S_j\setminus\{p\})$, $j=1,2$, in
Section~\ref{sec:geometry}. Set
\[
 s(z)=z+\chi(\arg z)\bigl(t(z)-z\bigr)\quad(z\ne0),\qquad s(0)=0.
\]
By Section~\ref{sec:geometry}, $t(z)-z=O(|z|^{1+\alpha})$ and
$Dt(z)-I=O(|z|^\alpha)$. Since $|D(\chi(\arg z))|\le\|\chi'\|_\infty/|z|$, we have 
\[
 s(z)-z=O(|z|^{1+\alpha}),\qquad Ds(z)-I=O(|z|^\alpha).
\]
Near the two real rays $s=z$, so extending $s$ by the identity below the
real axis gives a $C^1$ map with $Ds(0)=I$. Thus $s$ is a local diffeomorphism, with $Ds$ and $D(s^{-1})$ uniformly bounded
after shrinking the neighborhood of $0$. Since the extension fixes the lower
half disk, the inverse image of a sufficiently small upper half disk with
its diameter is a neighborhood of $0$ in the closed upper half plane.

\begin{figure}[!htbp]
  \centering
  \includegraphics[width=.95\linewidth]{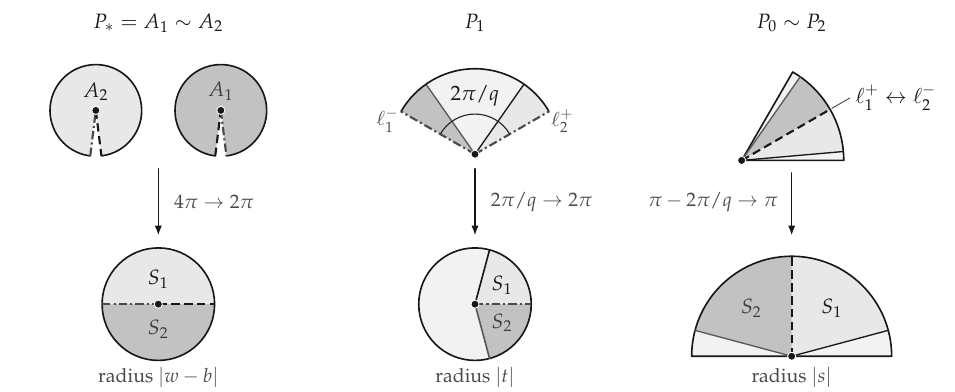}
  \caption{The angle changes at the three vertices of $M$, drawn for $q=3$.
  In each column the radial coordinate written below the chart is
  preserved. At $P_*$ the two slit disks are joined by
  \eqref{eq:surgery-bank-pairing}. In the upper right the wedge at $P_0$
  has been rotated by $-2\pi/q$ and joined to the wedge at $P_2$
  along $\ell_1^+\leftrightarrow\ell_2^-$.}
  \label{fig:surgery-vertex-charts}
\end{figure}

Along $\ell_2^-$ and $\ell_1^+$ the $t$-arguments are $\varphi$ and
$\varphi+2\pi/q$, respectively. Since $t(z)/z=1+O(|z|^\alpha)$, the
$z$-arguments tend to these same angles as $z\to0$. For $0<|z|$ sufficiently small,
$\chi(\arg z)=1$ on a neighborhood of each bank, so $s=t$ there. In the
$s$-plane the sector at $P_2$ has angles $[0,\varphi]$, and the sector at $P_0$ has angles
$[\varphi+2\pi/q,\pi]$. By \eqref{eq:surgery-bank-pairing}, paired points on
$\ell_2^-$ and $\ell_1^+$ have the same $w=t^q$ and hence the same radius
$|s|=|w|^{1/q}$. Rotating the sector at $P_0$ by $-2\pi/q$ thus identifies
the paired points on the ray of angle $\varphi$ and gives a wedge of opening $\pi-2\pi/q$. Put $a=\pi/(\pi-2\pi/q)$. Send $P_0\sim P_2$ to $0$ and use the coordinate
\[
 s\longmapsto
 \begin{cases}
  |s|e^{\sqrt{-1}a\arg s}, & \text{in the sector at }P_2,\\
  |s|e^{\sqrt{-1}a(\arg s-2\pi/q)}, & \text{in the sector at }P_0.
 \end{cases}
\]
At each fixed positive radius, the image angles run from $0$ to $a\varphi$
in the sector at $P_2$ and from $a\varphi$ to $\pi$ in the sector at $P_0$.
The two formulas agree on the paired banks and map the original boundary
rays onto the two halves of the real diameter. Since the radius $|s|$ is
preserved, this gives a homeomorphism from a neighborhood of $P_0\sim P_2$
in $M$ onto an upper half disk.
Figure~\ref{fig:surgery-vertex-charts} shows the angle changes at $P_*$, $P_1$, and $P_0\sim P_2$.

At all three vertices the charts extend homeomorphically by assigning
zero to the vertex. On each open piece, the differential of the map
$(r,\vartheta)\mapsto(r,a\vartheta)$ has matrix $\operatorname{diag}(1,a)$
in the orthonormal polar frames at the source and target points. Here the factors are $1/2$,
$q$, and $\pi/(\pi-2\pi/q)$, respectively. Thus the derivatives of the charts and their inverses are bounded on each cut piece. Away from these three vertices, we extend the angle changing formulas smoothly across the paired banks
in the rotated $t$-charts. Near $P_0\sim P_2$, the two arcs of $\partial M$ inherited from $\partial\D\setminus\{p\}$ correspond in the $z$-coordinate to $z=x<0$ on the $P_0$ side and $z=x>0$ on the $P_2$ side, respectively. Near each arc, $s=z$, so the polar formula for the chart at $P_0\sim P_2$ is a smooth change of coordinates up to the boundary wherever $x\ne0$. These charts give the topological disk $M$ a smooth structure with boundary. Choose an orientation preserving diffeomorphism $\kappa\colon\clD\to M$. On a finite cover by relatively compact subcharts, the coordinate
expressions of $\kappa$ and $\kappa^{-1}$ have bounded derivatives. Combined with the estimates for these chart, this gives bounded derivatives and inverse derivatives for the changes from the original cut pieces to
the disk coordinates given by $\kappa$.

\begin{figure}[!htbp]
  \centering
  \includegraphics[width=.95\linewidth]{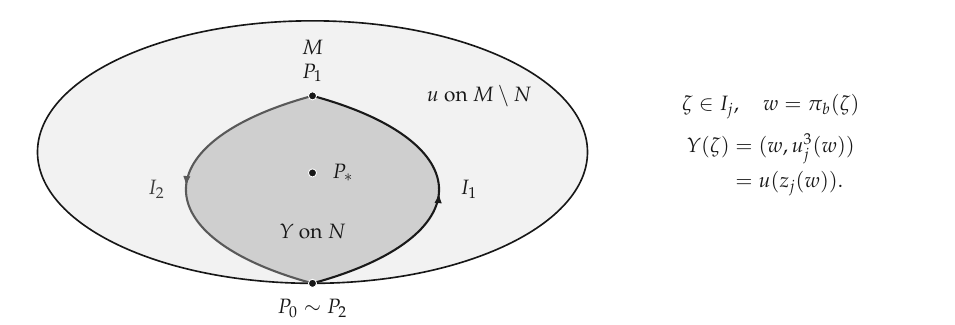}
  \caption{The reglued sectors form the disk $N\subset M$.
  The replacement equals $Y$ on $N$ and is induced by $u$ on
  $M\setminus N$. The two traces agree on $I_1$ and $I_2$.
  Here $P_1\in\partial N\cap\operatorname{int}M$ and
  $N\cap\partial M=\{P_0\sim P_2\}$.}
  \label{fig:surgery-replacement}
\end{figure}

Let $N\subset M$ be the closed disk obtained by regluing the cut sectors
$S_1$ and $S_2$ according to \eqref{eq:surgery-bank-pairing}. For
$x\in N$, let $\zeta(x)\in\clD$ be its coordinate in the covering disk,
so that $w=\pi_b(\zeta(x))$ is its projection to $\overline\Omega$.
Define the replacement on $N$ by $x\mapsto Y(\zeta(x))$. On $\overline{M\setminus N}$ we use the original map $u$. To specify
its argument more precisely, consider the projection from the cut disk to the original
parameter disk $\clD$ that sends each copy of a point to that point.
In the closure of the complement of the two cut sectors, the only distinct
points identified by \eqref{eq:surgery-bank-pairing} are $P_0$ and $P_2$,
and both project to $p$. The projection therefore induces a continuous
map $\rho\colon\overline{M\setminus N}\to\clD$. Define the replacement
on this closed set by $x\mapsto u(\rho(x))$. The two closed sets meet along the boundary of the covering disk $N$,
which in the $\zeta$-coordinate consists of $I_1$, $I_2$, and their two
common endpoints. If $\zeta(x)\in I_j$ and $w=\pi_b(\zeta(x))$, then
$\rho(x)$ is the point of $\partial S_j\setminus\{p\}$ with local
coordinate $z_j(w)$. By the definition of the graph height $u^3_j$ and
the choice of boundary values for $Y$ on $I_j$, we have
\[
 Y(\zeta(x))=(w,u^3_j(w))=u(\rho(x)),\qquad j=1,2.
\]
For $x=P_1$ and for $x=P_0\sim P_2$, we have $\rho(x)=p$ and
$\pi_b(\zeta(x))=0$. Since $u^3_j(0)=0$ for $j=1,2$, both maps take
the value $u(p)=0$ there. The two continuous maps thus agree on their
intersection and define a continuous map on $M$. Composing this map with
$\kappa\colon\clD\to M$ gives the comparison disk
$\widetilde u\colon\clD\to\R^3$.
Figure~\ref{fig:surgery-replacement} shows where the two maps are joined in $M$.

For $x\in\D\setminus(\ell_1\cup\ell_2)$, let $[x]\in M$ be the image of $x$
after cutting along the slits and identifying the banks by
\eqref{eq:surgery-bank-pairing}. The coordinate change from the original
parameter disk to the disk parametrizing $M$ is
\[
 F\colon\D\setminus(\ell_1\cup\ell_2)\longrightarrow\D,
 \qquad F(x)=\kappa^{-1}([x]).
\]
The restriction of $F$ to any open set $V\subset\D\setminus(\ell_1\cup\ell_2)$
is a diffeomorphism onto $F(V)$. A map $X$ on $V$ is therefore expressed
in disk coordinates as $X\circ F^{-1}$. Let $\sigma_1^2\ge\sigma_2^2>0$
be the eigenvalues of $(DF)^{\mathsf T}DF$ at $x$. The chart bounds above for $\kappa$ give
$\|(DF)^{-1}\|_{\mathrm{op}}^2 J_F=\sigma_1/\sigma_2\le C$ uniformly in $x$.
For $X\in W^{1,2}(V,\R^3)$ we therefore have
\begin{equation}\label{eq:surgery-coordinate-energy}
 E(X\circ F^{-1},F(V))
 =\frac12\int_V |DX(x)(DF(x))^{-1}|^2 J_F(x)\,dx
 \le C E(X,V).
\end{equation}
On $V=\D\setminus(S_1\cup S_2)$ we have $\rho([x])=x$, so
$\widetilde u=u\circ F^{-1}$ on $F(V)$. Taking $X=u|_V$ in
\eqref{eq:surgery-coordinate-energy} gives $E(\widetilde u,F(V))\le C E(u,V)<\infty$.

For the two sectors forming $N$, we first compute the energy in the
$w$-coordinate. Over $\Omega\setminus[0,b]$ the cover $\pi_b$ has two
sheets. If $u^3_-\not\equiv0$, choose one inverse branch $\sigma$ of
$\pi_b$ there and set $v=g\circ\widetilde T\circ\sigma$ as in
\eqref{eq:energytransport}. By \eqref{eq:newmap}, the two expressions for
$Y$ are $w\mapsto(w,u^3_+(w)+v(w))$ and $w\mapsto(w,u^3_+(w)-v(w))$.
If $u^3_-\equiv0$, the same expressions hold with $v=0$. Thus, in each case, the sum of the two
energies in the $w$-coordinate is
\begin{align*}
 \frac12\int_\Omega
 \bigl(4+|\nabla(u^3_++v)|^2+|\nabla(u^3_+-v)|^2\bigr)\,dw =2|\Omega|+\int_\Omega
 \bigl(|\nabla u^3_+|^2+|\nabla v|^2\bigr)\,dw.
\end{align*}
Exchanging the two sheets replaces $v$ by $-v$ and preserves the sum.
By \eqref{eq:decay}, $|\nabla u^3_+(w)|\le C|w|^{\alpha/q}$. When
$u^3_-\not\equiv0$, \eqref{eq:maj-ellipticity} and \eqref{eq:energytransport}
give $\int_\Omega |\nabla v|^2\,dw\le 2E(g)<\infty$. Thus the sum of two energies is
finite in both cases. On each open slit sector, $t\mapsto w=t^q$ is conformal
onto $\Omega\setminus[0,b]$, so the energy in $t$ equals the energy in $w$.
Since $z\mapsto t(z)$ is bi-Lipschitz by Section~\ref{sec:geometry}, the replacement on $N$ has finite energy on each slit sector in the original parameter
disk. For each of the two open slit sectors $V=\operatorname{int}(S_1)\setminus\ell_2$
and $V=\operatorname{int}(S_2)\setminus\ell_1$,
\eqref{eq:surgery-coordinate-energy} gives $E(\widetilde u,F(V))<\infty$.

Let $x\in\partial N$ with $\zeta(x)\in I_j$, and suppose that
$w=\pi_b(\zeta(x))$ lies on an open radial side or on the open circular
arc of $\partial\Omega$. The traces of $\widetilde u$ at $\kappa^{-1}(x)$ from
$\kappa^{-1}(\operatorname{int}N)$ and $\kappa^{-1}(M\setminus N)$ are
$Y(\zeta(x))$ and $u(\rho(x))$, respectively. Both equal $(w,u^3_j(w))$. Next let $x\in N$ be
the image of two bank points paired in \eqref{eq:surgery-bank-pairing}
over $w=\lambda b$, $0<\lambda<1$. The traces at $\kappa^{-1}(x)$ from
the two slit sectors both equal $Y(\zeta(x))$ by continuity of $Y$ on the
covering disk. In either case the interface divides a sufficiently small
disk about $\kappa^{-1}(x)$ into two Lipschitz domains. For a test
function supported in the disk, the boundary terms from integration by
parts on the two domains cancel along their common arc. Thus the piecewise derivative of
$\widetilde u$ is its weak derivative on the disk. Let $\Sigma\subset\operatorname{int}M$ consist of $P_*$, $P_1$, and the
four points $x\in\partial N$ for which
$$
\pi_b(\zeta(x))\in\{r e^{\sqrt{-1}\theta_0},r e^{\sqrt{-1}(\theta_0+\theta)}\}.
$$
Each of these two corners has one preimage on each $I_j$. Of the two points
over $0$, $P_1$ lies in $\operatorname{int}M$ and $P_0\sim P_2$ lies on
$\partial M$. We have now shown
that $\widetilde u\in W^{1,2}_{\mathrm{loc}}(\D\setminus\kappa^{-1}(\Sigma),\R^3)$.
Since the interfaces have planar measure zero, summing the three energy
bounds obtained from \eqref{eq:surgery-coordinate-energy} gives
$\int_{\D\setminus\kappa^{-1}(\Sigma)}|D\widetilde u|^2<\infty$.

For $a\in\kappa^{-1}(\Sigma)$, choose a disk $B$ centered at $a$ with
$\overline B\subset\D$ and $B\cap\kappa^{-1}(\Sigma)=\{a\}$. For $\phi\in C_c^\infty(B)$, test
the weak derivative identity on $B\setminus\{a\}$ with
$\eta_\eps(\cdot-a)\phi$, where $\eta_\eps$ is given by
\eqref{eq:transfer-log-cutoff}. For $k=1,2$, the term containing the
derivative of the cutoff satisfies
\[
 \left|\int_B \widetilde u(y)\phi(y)\,\partial_k\eta_\eps(y-a)\,dy\right|
 \le \|\widetilde u\|_{L^\infty(B)}\|\phi\|_{L^2(B)}
       \left(\frac{2\pi}{\log(1/\eps)}\right)^{1/2}\longrightarrow0.
\]
Since $0\le\eta_\eps\le1$ and $\eta_\eps(y-a)\to1$ for $y\ne a$, letting
$\eps\searrow0$ gives $\int_B\widetilde u\,\partial_k\phi
=-\int_B(\partial_k\widetilde u)\phi$. Hence the $L^2$ integrable piecewise derivative is the weak derivative
of $\widetilde u$ throughout $\D$. Recalling that $\widetilde u$ is continuous
on $\clD$, we conclude that
$\widetilde u\in C^0(\clD,\R^3)\cap W^{1,2}(\D,\R^3)$.

In the positive boundary traversal in \eqref{eq:surgery-bank-order}, the
arc inherited from $\partial\D\setminus\{p\}$ runs from $P_2$ to $P_0$
in the original boundary orientation. After the bank identifications in
\eqref{eq:surgery-bank-pairing}, this arc with its endpoints joined is
$\partial M$. Thus $\rho|_{\partial M}\colon\partial M\to\partial\D$ is
an orientation-preserving homeomorphism taking $P_0\sim P_2$ to $p$.
Since $\kappa$ also preserves orientation, the trace
\[
 \widetilde u|_{\partial\D}
 =(u|_{\partial\D})\circ(\rho|_{\partial M})\circ(\kappa|_{\partial\D})
\]
is a weakly monotone parametrization of $\Gamma$. Hence
$\widetilde u\in\A(\Gamma)$.

On $V=\D\setminus(S_1\cup S_2)$, the identity $\widetilde u\circ F=u$
gives $A(\widetilde u,F(V))=A(u,V)$. For each $j=1,2$, changing
variables by $z_j$ and using $u(z_j(w))=(w,u^3_j(w))$ gives
\[
 A(u,\operatorname{int}S_j)=\int_\Omega W(\nabla u^3_j)\,dw=A(u^3_j).
\]
On the two open slit sectors in $\kappa^{-1}(\operatorname{int}N)$,
we have $\widetilde u=Y\circ\zeta\circ\kappa$. The coordinate map
$\zeta\circ\kappa$ takes these sectors diffeomorphically onto the two
sheets of $\pi_b$ over $\Omega\setminus[0,b]$. Since the two seams over $[0,b]$ and the boundary of $N$ have planar
measure zero in disk coordinates, these two sectors together contribute $A(Y)$
to $A(\widetilde u)$. Adding the areas gives
\begin{equation}
\label{eq:sewn}
 A(\widetilde u)=A(u)-A(u^3_1)-A(u^3_2)+A(Y).
\end{equation}
For the minimizer $u$ in Theorem~\ref{thm:main}, the admissibility of
$\widetilde u$ and \eqref{eq:mine} imply $A(u)\le A(\widetilde u)$.
If $u^3_-\not\equiv0$ on $\Omega$, however, \eqref{eq:sewn} and
\eqref{eq:drop} give
\[
 A(\widetilde u)-A(u)=A(Y)-A(u^3_1)-A(u^3_2)<0.
\]
This contradiction rules out the assumed boundary branch point $p$
whenever the two graph heights differ somewhere on $\Omega$.

\subsection{The coincident graphs}\label{sec:surgery-coincident}

In this subsection, we exclude the remaining case $u^3_-\equiv0$ on
$\Omega$. The replacement $\widetilde u$ constructed in
Subsection~\ref{sec:surgery-sobolev} then preserves area. We show that a
conformal reparametrization of $\widetilde u$ minimizes $E$ but fails to be locally
injective at the interior point corresponding to $P_*$, contrary to
the interior immersion theorem in~\cite[Theorem~6.2]{gull}. This completes
the proof of Theorem~\ref{thm:main}.

Since $u^3_1=u^3_2=u^3_+$, the map $Y=(\pi_b,u^3_+\circ\pi_b)$
parametrizes the same graph on both sheets. Thus
$A(Y)=A(u^3_1)+A(u^3_2)$ by Subsection~\ref{sec:transfer}, and
\eqref{eq:sewn} gives $A(\widetilde u)=A(u)$. Moreover,
$\widetilde u\circ F=u$ on both open slit sectors, since on each sector
the point over $w$ has value $(w,u^3_+(w))=u(z_j(w))$. The same identity
holds on $\D\setminus(S_1\cup S_2)$ by construction of $\widetilde u$ and $F$. By
Subsection~\ref{sec:surgery-sobolev}, the weak derivative of $\widetilde u$
agrees with its derivative on the images under $F$ of these three open
sets. In particular,
\[
 D\widetilde u(F(x))=Du(x)(DF(x))^{-1}
 \qquad\text{for a.e. }x\in\D.
\]
Both factors on the right are uniformly bounded by the $C^1$ regularity
of $u$ and the above estimates for $F$. Hence $\widetilde u\in
W^{1,\infty}(\D,\R^3)$, and its continuous representative is Lipschitz on
$\clD$.

To choose conformal parameters on $\D$, set $G=(D\widetilde u)^{\mathsf T}D\widetilde u$.
By conformality of $u$, $Du(x)^{\mathsf T}Du(x)=\tfrac12|Du(x)|^2 I$.
Thus $G(F(x))$ vanishes where $Du(x)=0$ and elsewhere is a positive
multiple of $\big((DF(x))^{-1}\big)^{\mathsf T}(DF(x))^{-1}$. In the latter case,
with $\sigma_1^2,\sigma_2^2$ the eigenvalues of $DF(x)^{\mathsf T}DF(x)$
used to obtain \eqref{eq:surgery-coordinate-energy}, we have
\[
 \sqrt{\frac{\lambda_{\max}(G(F(x)))}{\lambda_{\min}(G(F(x)))}}
 =\frac{\sigma_1}{\sigma_2}.
\]
The previous bound for $\sigma_1/\sigma_2$ given by $\kappa$ therefore makes the conductivity
$\sigma_G=\sqrt{\det G}\,G^{-1}$ uniformly elliptic on the set where $G$ is
positive definite. Take $\sigma_G=I$ elsewhere and extend it by $I$ off $\D$. Since
$\det\sigma_G=1$, applying the conductivity straightening result in~\cite[Lemma~3.1]{alp}
and postcomposing with a Riemann map as in Subsection~\ref{sec:transfer}
gives an orientation preserving quasiconformal map $\Psi\colon\D\to\D$
that extends to a homeomorphism of $\clD$ and satisfies
\[
 D\Psi\,\sigma_G\,(D\Psi)^{\mathsf T}=J_\Psi I\qquad\text{a.e. in }\D.
\]

Set $u_*=\widetilde u\circ\Psi^{-1}$. Since $\widetilde u$ is Lipschitz
and $\Psi^{-1}\in W^{1,2}(\D,\R^2)$, we have $u_*\in W^{1,2}(\D,\R^3)$.
It is continuous on $\clD$, and its trace is an orientation preserving reparametrization of $u|_{\partial\D}$. In particular, $u_*\in\A(\Gamma)$.
At almost every $y$ with $G(y)$ positive definite, the above identity for $D\Psi$ gives
\[
 Du_*(\Psi(y))^{\mathsf T}Du_*(\Psi(y))
 = ((D\Psi(y))^{-1})^{\mathsf T}G(y)(D\Psi(y))^{-1}
 = \frac{\sqrt{\det G(y)}}{J_\Psi(y)}I.
\]
Where $G(y)=0$, we have $Du_*(\Psi(y))=0$ almost everywhere as well.
Thus $u_*$ is weakly conformal and changing variables by $\Psi$ in the
area  functional gives $A(u_*)=A(\widetilde u)$.
Together with \eqref{eq:mine}, this yields
\[
 E(u_*)=A(u_*)=A(\widetilde u)=A(u)=\inf_{v\in\A(\Gamma)}E(v).
\]

We apply the interior immersion theorem in~\cite[Theorem~6.2]{gull} with
$H=0$ to both $u$ and $u_*$. Since disk reflection preserves the Dirichlet
energy and exchanges the two boundary orientations, both maps minimize among disks with monotone trace onto $\Gamma$ in either orientation and have homeomorphic boundary traces and are harmonic in $\D$. For each
map, some component of its complex derivative is holomorphic and not
identically zero, since its trace parametrizes $\Gamma$. Its Dirichlet
energy density is therefore positive except possibly at isolated points.
It follows that $u$ and $u_*$ are immersions in $\D$.

\begin{figure}[!htbp]
  \centering
  \includegraphics[width=.95\linewidth]{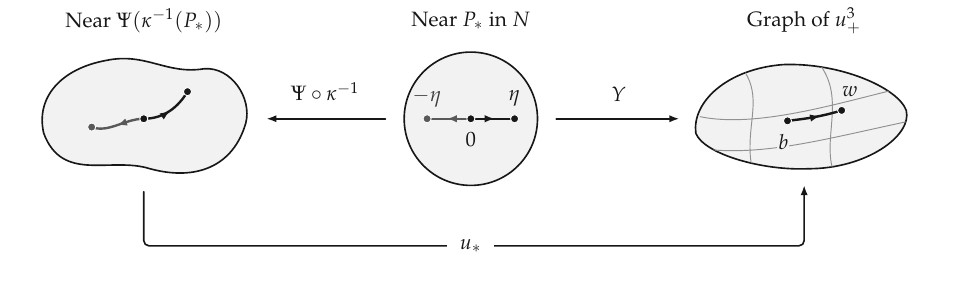}
  \caption{The common graph when $u^3_-\equiv0$. Near $P_*$ the points
  with coordinates $\eta$ and $-\eta$ have the same image $(w,u^3_+(w))$
  under $Y$, where $w=b+\eta^2$. The right-hand labels give horizontal
  coordinates. For $\eta\ne0$, the homeomorphism $\Psi\circ\kappa^{-1}$
  gives two distinct points with the same image under $u_*$. Both
  approach $\Psi(\kappa^{-1}(P_*))$ as $\eta\to0$.}
  \label{fig:coincident-graph}
\end{figure}

Consider now the interior point $\Psi(\kappa^{-1}(P_*))$ in the domain
of $u_*$. In the coordinate $\eta$ centered at $P_*$, the projection
is $w=b+\eta^2$. For every sufficiently small $\eta\ne0$, the distinct
points of $N$ with coordinates $\eta$ and $-\eta$ both project to $w$
and hence both have value $(w,u^3_+(w))$ under the replacement. Their
images under $\Psi\circ\kappa^{-1}$ are distinct points with the same
$u_*$-value, and both tend to $\Psi(\kappa^{-1}(P_*))$ as $\eta\to0$,
as shown in Figure~\ref{fig:coincident-graph}. Thus $u_*$ fails to be
locally injective at an interior point, contradicting the immersion
conclusion. This excludes the assumed boundary branch point $p$ when
$u^3_-\equiv0$.

\begin{proof}[Proof of Theorem~\ref{thm:main}]
Let $u$ minimize $E$ over $\A(\Gamma)$. By \eqref{eq:mine}, $u$ is a
least-area disk in $\R^3$. Such a disk has neither true interior branch
points~\cite{osserman} nor false interior branch points~\cite{gull,alt1,alt2}.
Thus $u$ is an immersion in $\D$.

Suppose $du(p)=0$ for some $p\in\partial\D$, and let $u^3_1,u^3_2$ be
the two graph heights over the sector $\Omega$ in Section~\ref{sec:geometry}.
If $u^3_1\not\equiv u^3_2$, then \eqref{eq:drop} and \eqref{eq:sewn} give
$A(\widetilde u)<A(u)$ for the admissible disk $\widetilde u$ constructed
in Subsection~\ref{sec:surgery-sobolev}, contrary to the minimality of $u$.
If $u^3_1\equiv u^3_2$, the disk $u_*$ constructed in
Subsection~\ref{sec:surgery-coincident} is also a least-area disk.
Near the interior point corresponding to $P_*$, it is a branched double
cover of the common graph and hence fails to be locally injective there.
This contradicts the above interior immersion result applied to $u_*$. Thus
$du(p)\ne0$ for every $p\in\partial\D$. By Section~\ref{sec:geometry}, $u\in C^{1,\alpha}(\clD,\R^3)$, so the
conformality equations of $u$ extend to $\partial\D$. Since $du$ is nonzero
there, it has rank two at every boundary point. Hence $u$ is a
$C^{1,\alpha}$ conformal immersion on $\clD$. In particular, the boundary homeomorphism $u|_{\partial\D}$ has nonvanishing tangential derivative and is a $C^{1,\alpha}$ diffeomorphism onto $\Gamma$.
\end{proof}

\begin{example}\label{ex:planar-c1}
We construct a planar least-area disk that has a boundary branch point when it spans a regular $C^1$ Jordan curve.

Fix $0<r<e^{-\pi}$. On $\overline{D_r^+}\setminus\{0\}$ take the branch
of $\log z$ with $0\le\arg z\le\pi$, and set
\[
 f(z)=-\frac{z}{\log z}\quad(z\ne0),\qquad f(0)=0.
\]
Since $f(z)=o(|z|)$ and $f'(z)=-1/\log z+1/(\log z)^2\to0$ as
$z\to0$, we have $f\in C^1(\overline{D_r^+})$ and $df(0)=0$.
For $z\ne0$, set $\vartheta=\arg z$. Since $-\log|z|>\pi\ge\vartheta$,
we have
\[
 \operatorname{Re}f'(z)=\frac{-\log|z|}{(\log|z|)^2+\vartheta^2}
       +\frac{(\log|z|)^2-\vartheta^2}{\bigl((\log|z|)^2+\vartheta^2\bigr)^2}>0.
\]
For distinct $z_1,z_2\in\overline{D_r^+}$,
\[
 \operatorname{Re}\frac{f(z_2)-f(z_1)}{z_2-z_1}
 =\int_0^1\operatorname{Re}f'\bigl(z_1+t(z_2-z_1)\bigr)\,dt>0.
\]
Thus $f$ is injective on the closed half disk.

Choose a smooth Jordan domain $U\subset D_r^+$ such that
$U\cap B_\eps(0)=D_\eps^+$ for some $0<\eps<r$. Identify $\C$ with
the coordinate plane in $\R^3$ and give $\Gamma_0=f(\partial U)$ the
positive orientation as the boundary of $f(U)$. The restriction of
$f$ to $\partial U$ is a homeomorphism onto $\Gamma_0$. Away from $0$,
$f$ is locally holomorphic with nonzero derivative, so this Jordan
curve is smooth there. On $(-\eps,\eps)\subset\partial U$, the coordinate
$s=\operatorname{Re}f(x)$ is strictly increasing. Write this part of
$\Gamma_0$ as the graph $(s,y(s))$, where $y(s)=\operatorname{Im}f(x)$.
For $x>0$ we have $y(s)=0$. For $x<0$, the choice $\arg x=\pi$ gives
\[
 f(x)=\frac{x(-\log|x|+\sqrt{-1}\pi)}{(\log|x|)^2+\pi^2},\qquad
 s=\frac{x(-\log|x|)}{(\log|x|)^2+\pi^2},\qquad
 y(s)=\frac{\pi s}{-\log|x|}.
\]
In particular, $y(s)/s\to0$ as $s\to0$, so $y'(0)=0$. Differentiating
the real and imaginary parts for $x<0$ gives
\[
 y'(s)=\frac{\operatorname{Im}f'(x)}{\operatorname{Re}f'(x)}
       \sim\frac{\pi}{-\log|x|}\qquad\text{as }x\nearrow0.
\]
Thus $s\mapsto(s,y(s))$ is a regular $C^1$ parametrization through the
origin, with derivative $(1,0)$ there.

For $s<0$, the arc length $t$ from the origin to $(s,y(s))$ satisfies
$t=\int_s^0\sqrt{1+y'(\sigma)^2}\,d\sigma\sim|s|$. Since
$|s|\sim |x|/(-\log|x|)$, we have $\log(1/t)\sim-\log|x|$. The angle between the
tangent at $(s,y(s))$ and the horizontal tangent at the origin is
$\arctan y'(s)\sim\pi/\log(1/t)$. The tangent fails the Dini condition at
$0$, since $\int_0^\delta dt/(t\log(1/t))=\infty$ for $0<\delta<1$.
Figure~\ref{fig:planar-c1-boundary} shows the boundary arc at two scales.

\begin{figure}[!htbp]
  \centering
  \includegraphics[width=.93\linewidth]{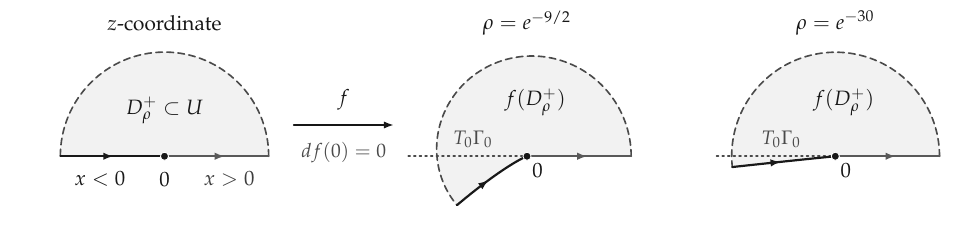}
  \caption{The real diameter $[-\rho,\rho]$ and its image under
  $f(z)=-z/\log z$ are drawn with solid lines. Dashed lines mark the
  upper semicircle of radius $\rho$ and its image. In the two image
  panels, $[0,f(\rho)]$ has the same displayed length. As $\rho$
  decreases, $f([-\rho,0])$ flattens toward the horizontal tangent
  $T_0\Gamma_0$.}
  \label{fig:planar-c1-boundary}
\end{figure}

Let $\varphi\colon\D\to U$ be a Riemann map. Since $\partial U$ is
smooth, the boundary regularity theorem for conformal maps in~\cite[Theorems~1 and~2]{warschawski}
implies that $\varphi$ extends to a $C^1$ diffeomorphism of $\clD$
onto $\overline U$.
The map $u_0=f\circ\varphi$ parametrizes $f(U)$ conformally and maps
$\partial\D$ homeomorphically onto $\Gamma_0$ with positive orientation.
Hence $u_0\in\A(\Gamma_0)$ and $E(u_0)=A(u_0)=|f(U)|$. At the unique
point $p_0=\varphi^{-1}(0)\in\partial\D$, its differential is
$du_0(p_0)=df(0)\,d\varphi(p_0)=0$.

Let $v\in\A(\Gamma_0)$, and let $w\colon\clD\to\C$ be the continuous
harmonic extension of $v|_{\partial\D}$. The orthogonal projection of
$v$ onto $\C$ has the same trace and energy at most $E(v)$. The
Dirichlet principle therefore gives $w\in W^{1,2}(\D,\C)$ and
$E(w)\le E(v)$. For each $a\in f(U)$, the weakly monotone trace
$v|_{\partial\D}$ traverses $\Gamma_0$ once with positive orientation,
so its winding number about $a$ is one. If $a\notin w(\D)$, then $w$
would extend the trace $v|_{\partial \D}$ to a disk in $\C\setminus\{a\}$, forcing its
winding number about $a$ to be zero. Thus $f(U)\subset w(\D)$. Since
$w$ is smooth in $\D$, the area formula counts every point of $f(U)$
with multiplicity at least one. For the corresponding surjectivity and
area bound in a planar metric, see Hildebrandt--von der Mosel~\cite[proof of Theorem~2, (7)--(10)]{hvm05}.
It follows that
\[
 E(v)\ge E(w)\ge A(w)\ge |f(U)|=E(u_0).
\]
Thus $u_0$ minimizes $E$ in $\A(\Gamma_0)$ and satisfies $du_0(p_0) = 0$. 
\end{example}

\vspace{0.5cm}

\bibliographystyle{amsalpha}
\bibliography{references}
\end{document}